\documentclass[12pt]{article}   
\usepackage{amsfonts}
\usepackage{latexsym, amsmath, amssymb, amsthm, mathrsfs}
\usepackage{color}
\usepackage{multirow}
\usepackage{graphicx}

\newcommand{\ai}{\'{\i}}   
\newcommand{\eps}{\epsilon}

\newcommand{\uno}{\mbox{\large 1}\hspace{-4pt}\mbox{\large l}}
\newcommand{\neno}{\newline\noindent}

\newcommand{\och}{\mbox{o}_{\Pr}}
\newcommand{\ip}{\mbox{\small(p)}}
\newcommand{\id}{\mbox{\small(d)}}
\newcommand{\cip}{\stackrel{\ip}{\to}}
\newcommand{\cid}{\stackrel{\id}{\to}}

\newcommand{\calc}{{\cal{C}}}

\newcommand{\calm}{{\cal{M}}}

\newcommand{\calk}{{\cal{K}}}
\newcommand{\cals}{{\cal{S}}}
\newcommand{\calh}{{\cal{H}}}
\newcommand{\calz}{{\cal{Z}}}

\newcommand{\real}{\mathbb{R}}

\newcommand{\beq}{\begin{equation}}   
\newcommand{\enq}{\end{equation}}

\newcommand{\os}{\overline{S}}
\newcommand{\hs}{\hat{S}}
\newcommand{\beqa}{\begin{eqnarray}}   
\newcommand{\enqa}{\end{eqnarray}}

\newcommand{\ma}{\mathbf{a}}
\newcommand{\mb}{\mathbf{b}}
\newcommand{\hamu}{\hat{\mu}}

\newtheorem{theo}{Theorem}

\newtheorem{lem}{Lemma}

\graphicspath{{./PruebasPotenciaEstadisticos/}}

\begin{document}
\begin{center}
\begin{Large}
	Using the Sinkhorn divergence in permutation tests for the multivariate two-sample problem. 
\end{Large}
\vskip 2pc
\begin{large}
E. Del Barrio,\footnote{Email: eustasio.delbarrio@uva.es} 
J. S. Osorio \footnote{Email: js.osorio125@uniandes.edu.co} and A. J. Quiroz \footnote{Email: ajquiroz@gmail.com}  
\end{large}
\end{center}
\vskip 1pc
\begin{abstract}
In order to adapt the Wasserstein distance to the large sample multivariate non-parametric 
two-sample problem, making its application computationally feasible, permutation 
tests based on the Sinkhorn divergence between probability vectors associated to 
data dependent partitions are considered. Different ways of implementing these tests 
are evaluated and the asymptotic distribution of the underlying statistic is established in some cases. 
The statistics proposed are compared, in simulated examples, with the test of
Schilling's, one of the best non-parametric tests available in the literature.

\end{abstract}
\pagebreak	
\section{Introduction}
The Wasserstein distance has been established as a powerful tool in diverse problems in 
the context of multivariate statistics. Del Barrio et al, \cite{Del Barrio}, discuss 
statistics for goodness of fit based on the $L^2$ distance between an empirical distribution 
and the hypothesized model. Frogner et al, \cite{frogner}, consider applications to 
regression, while Sommerfeld and Munk, \cite{sommerfeld}, describe applications to 
supervised learning of the Wasserstein distance. References to several other applications 
of the Wasserstein distance to diverse data mining problems are given by Mena and Weed, 
\cite{mena}. The Wasserstein distance is the solution of the (classical or unregularized) 
Optimal Transport (OT) problem, as described in Chapter 1 of Villani's book, \cite{villani}, 
and we will refer to it using these terms as well.
 
In the non-parametric multivariate two-sample problem, we have samples 
$X_1,\dots,X_m$ i.i.d. from a distribution $P$ and  $Y_1,\dots,Y_n$ i.i.d. from a 
distribution $Q$, both distributions defined on $\real^d$. The null hypothesis to be tested 
is $P=Q$, against the alternative $P\neq Q$. There are a few true non-parametric tests 
available in the literature for this problem, with the main ones, in terms of power 
against the alternative, being of a graph theoretic nature. Friedman and Rafsky, \cite{fr79}, 
propose several ground breaking graph theoretic options for this problem, including one based 
on the Minimal Spanning Tree (or multiple orthogonal spanning trees) for the joint sample 
(the union of the $X$ and $Y$ samples). Schilling, \cite{schilling}, studies a very powerful 
test based on the $k$ nearest neighbor graph of the joint sample and, more recently, Chen and 
Friedman, \cite{Chen}, consider a new kind of graph theoretic test for the non-parametric 
two-sample problem. Interestingly, all of these graph theoretic tests can be presented as 
permutation tests and that viewpoint is useful for working out their asymptotic 
distributions.  
 
In view of its success in other important statistical problems, it is natural
to try to use the Wasserstein distance in the two-sample problem. Still, direct application 
of OT in the setting of the multivariate two-sample problem is hampered by computational cost 
in the case of medium sized or large data sets. For instance, at the present time, computation 
of the Wasserstein distance between two samples of size 10,000 in dimension 5, takes about 70 
minutes on a laptop computer using an INTEL$^\copyright$ Core I5 processor, and for a single 
test, this calculation would have to be performed several times, since the null distribution 
of the statistic depends on the underlying common distribution, which is assumed unknown. 
The obstacle here is the size of the optimization problem that must be solved to obtain the OT 
between two large empirical distributions.

The Sinkhorn divergence is a modification of the Wasserstein metric, that offers important 
advantages in terms of computational speed and, as established recently by Geneway et al., 
\cite{genevay19}, also with respect to speed of convergence as a function of sample size, since 
in any dimension the Sinkhorn divergence between an empirical measure and its population 
counterpart decreases at the rate $1/\sqrt{n}$, avoiding the ``curse of dimensionality''. In 
this regard, see also Mena and Weed, \cite{mena}, who improve on the results of \cite{genevay19}, 
following the method of analysis developed by Del Barrio and Loubes, \cite{Del Barrio L}, for the 
unregularized OT problem.

The Sinkhorn divergence is obtained by adding an entropy cost term to the objective function 
that must be minimized in the optimal transport problem (details are given in Section 2). This 
modification comes at the cost of losing the metric properties enjoyed by the Wasserstein distance. The Sinkhorn divergence does not satisfy the axioms of a distance and it can occur that 
it takes a negative value. 

Bigot et al., \cite{bigot}, consider the asymptotic distribution of the Sinkhorn divergence 
between probability distributions on finite spaces and, as an example of a statistical 
application, they propose to use the Sinkhorn divergence in a bootstrap test for the two-sample 
problem. In their test, the Sinkhorn divergence is calculated between frequency (probability) 
vectors that, for each of the two samples, are associated to a fixed partition of the sample space. 
The key idea of this proposal is crucial. The main computational savings come from the fact that 
the original samples, of size $m$ and $n$, are replaced by frequency vectors of the same size as 
the partition chosen. 

In the present article, we consider modifying the procedure of Bigot et al., \cite{bigot}, 
described above, in two directions: (i) Use a permutation test instead of 
a bootstrap test and (ii) Use sample dependent partitions, such as those determined by $k$-means 
clustering of the joint sample, to partition the data in a way that, under the null hypothesis, 
converges, as the sample size grows, to a limiting ``natural'' partition of the sample space. Good, 
\cite{good05}, provides several reasons why permutations tests should be preferred over bootstrap 
tests. One such reason is that under the null hypothesis, a permutation test for the two sample 
problem is guaranteed to satisfy the nominal level of the test for finite
sample size (the only limitation being the size 
of the simulation carried out, that is, the number of random permutations generated). Current implementations of the $k$-means clustering procedure have computational complexity which is 
linear on the dimension, the number of clusters and the sample size, i.e. $O(k\,d(m+n))$ when applied to the joint sample (\cite{rokach05}). This linearity makes $k$-means very convenient for use on large data sets. General results of Pollard, \cite{pollard82}, on the 
convergence of the cluster centers for $k$-means clustering, together with tools from empirical processes theory, allow us to prove, under general conditions, that the distribution of the Sinkhorn 
divergence for our partition scheme under the null hypothesis, converges to a Gaussian limit 
depending only on the underlying common distribution and the limiting cluster centers. Other possible natural data dependent partition methods are considered
elsewhere, for reasons of space.

The rest of the article is organized as follows. Section 2 describes the Sinkhorn divergence (mostly in the case of finite spaces) and the different options of implementation of a two-sample 
statistic considered in our experiments. In Section 3, a central limit theorem is proved for the Sinkhorn divergence in our 
setting of data dependent partitions, while Section 4 describes the distributions (null and alternatives) 
considered in the evaluation experiments and discusses the results of those experiments.

\section{Options in permutation statistics based on Sinkhorn's divergence}
\subsection{Definitions of Wasserstein's distance and Sinkhorn's divergence}
For concreteness, we restrict the presentation of Wasserstein's distance and Sinkhorn's divergence 
to Euclidean space, with the usual norm, $\|\cdot\|$, although the discussion can be carried out on 
an arbitrary complete metric space. Let $\mu$ and $\nu$ be probability measures on $\real^d$. The 
($L^2$) Wasserstein (or Optimal Transport) distance between $\mu$ and $\nu$, $W_2(\mu,\nu)$, is 
given by 

\beq\label{wass}
W_2(\mu,\nu)=\inf_{\pi\in\Pi(\mu,\nu)}\left(\int_{\real^d\times\real^d}\|x-y\|^2\mbox{d}\pi(x,y)\right)^{1/2},
\enq 
where $\Pi(\mu,\nu)$ denotes the set of probability measures on $\real^d\times\real^d$ with margin\-als $\mu$ and $\nu$.
When $\mu$ and $\nu$ are discrete probability measures with finite support 
$\Omega=\{Z_1,\dots,Z_L\}\subset\real^d$, the probability measures can be identified with the
corresponding vectors of probabilities $\ma=(a_1,\dots, a_L)'$ and 
$\mb=(b_1,\dots, b_L)'$, that the measures $\mu$ and $\nu$, respectively, assign to 
the points in $\Omega$. In this context, the formula in (\ref{wass}) reduces to 
\beq\label{wass_discrete}
W_2^2(\ma,\mb)\>\>=\>\>\min_{T\in U(\ma,\mb)}\sum_{i,j\leq L}\|Z_i-Z_j\|^2T_{i,j}\>\>=
\>\>\min_{T\in U(\ma,\mb)}\left\langle D,T\right\rangle,
\enq 
where $T=(T_{i,j})_{i,j\leq L}$ is a $L\times L$ matrix of probabilities 
such \linebreak that $\sum_{i,j\leq L}T_{i,j}=1$, defining a probability distribution on
$\Omega\times\Omega$, $U(\ma,\mb)$ is the set of
probability matrices with marginals $\ma$ and $\mb$:
\linebreak $U(\ma,\mb)=\{T\in [0,1]^{L\times L}: T\uno=\ma, T'\uno=\mb\}$ with 
$\uno$ the vector of ones of length $L$, 
$D=(\|Z_i-Z_j\|^2)_{i,j\leq L}$ is the matrix of squared distances and 
$\left\langle \cdot,\cdot\right\rangle$ denotes the Frobenius inner product for
matrices. Sommerfeld and Munk, \cite{sommerfeld}, provide the asymptotic theory
for $W_2$ in (\ref{wass_discrete}) when $\mb=(a_{1,n},\dots,a_{K,n})$ is the empirical version of the probability vector $\ma$.

The entropic regularization of problem (\ref{wass_discrete}) is defined as follows:\linebreak
For a $L\times L$ probability matrix $T$, its Shannon entropy $H(T)$, is defined as \linebreak $H(T)=\sum_{i,j\leq L}T_{i,j}\ln(1/T_{i,j})$, in the understanding that when
$T_{i,j}=0$, the corresponding term in the sum is zero. For a cost parameter
$\lambda>0$, the entropy regularized OT problem, called Sinkhorn divergence, is
given by 

\beq\label{sinkhorn}
S_\lambda(\ma,\mb)\>=\>\min_{T\in U(\ma,\mb)}\left\langle D,T\right\rangle
-\lambda H(T),
\enq 
with $U(\ma,\mb)$ as defined before. $S_\lambda(\ma,\mb)$, as defined in (\ref{sinkhorn})
does not satisfy the axioms of a distance. In particular, $S_\lambda(\ma,\mb)$ can take negative values, specially for larger values of $\lambda$. 
To alleviate this problem, the following 
modifications have been proposed. The definition of Sinkhorn's divergence considered in \cite{cuturi} is
\beq\label{sinkhorn_3}
\hs_\lambda(\ma,\mb)=\left\langle D,T^\lambda\right\rangle,\>\mbox{ where }
\> T^\lambda=\mbox{argmin}_{T\in U(\ma,\mb)}\left\langle D,T\right\rangle
-\lambda H(T).
\enq
$\hs_\lambda(\ma,\mb)$ is always non-negative, but can fail to satisfy the 
coincidence axiom. On the other hand,
\beq\label{sinkhorn_2}
\os_\lambda(\ma,\mb)=\hs_\lambda(\ma,\mb)-\frac{1}{2}\left(\hs_\lambda(\ma,\ma)+\hs_\lambda(\mb,\mb)\right)
\enq
satisfies the coincidence axiom (see \cite{genevay19}) and is always non-negative.
In our evaluations, $S_\lambda(\ma,\mb)$, 
$\os_\lambda(\ma,\mb)$ and $\hs_\lambda(\ma,\mb)$ will be considered, together with the
limiting case of classical OT, obtained from $S_\lambda(\ma,\mb)$ by making $\lambda=0$.

There are three reasons that serve as motivation for considering the entropic 
regularization of the optimal transport problem (\ref{wass_discrete}). \begin{itemize}
\item[(i)] Computational speed. The solution of the classical OT problem 
has complexity $\mbox{O}(L^3\ln L)$, while the computational complexity of
the corresponding regularized problem is $\mbox{O}(L^2)$ (see \cite{cuturi}).
\item[(ii)] Smoothness of the solution. The solution of (\ref{wass_discrete})
occurs at a matrix $T^0$, whose entries are almost all zeroes (at most
$2L-1$ elements of $T^0$ can be different from $0$). The solution
to the regularized problem can be reached at a less extreme probability matrix
$T^\lambda$ (see the discussion in \cite{cuturi}).
\item[(iii)] Avoiding the ``curse of dimensionality''. As mentioned in the Introduction, Geneway et al., \cite{genevay19}, have proved that the the Sinkhorn divergence between an empirical measure and its population counterpart decreases at the rate $1/\sqrt{n}$ (with $n$ being the sample size) regardless of the dimension, 
while the rate of convergence of the classical OT distance is $\mbox{O}(1/n^{1/d})$, (see Theorem 11 in \cite{NilesWeedRigollet}.)
\end{itemize} 
The relevance of arguments (i) and (iii) in our particular application will be evaluated through the results of our experiments.

\subsection{Sinkhorn's divergence between frequency vectors determined by data
dependent partitions}
For $m,n\in\mathbb{N^+}$ and $N=n+m$, suppose i.i.d. samples $X_1,\dots,X_n$ and 
$Y_1,\dots,Y_m$ in $\real^d$ are available from 
the continuous distributions $P$ and $Q$, respectively. We will call these samples the $X$ and $Y$ samples. The joint sample, $Z_1,\dots, Z_N$ is obtained 
by concatenation of the $X$ and $Y$ samples. Assume both distributions 
have support contained in a set $\cals\subset\real^d$ and a partition
$\calc=\{C_1,\dots,C_k\}$ of $\cals$ is available. Let $\ma_N$ and $\mb_N$ denote, 
respectively, the vectors of frequencies for the $X$ and $Y$ samples respect to the 
sets of $\calc$. That is,  $\ma_N=(a_{1,N},\dots,a_{k,N})'$, where
\beq\label{freq}
a_{j,N}=\frac{\#\{i\leq n:X_i\in C_j\}}{n}, \>\mbox{ for each }\> j\leq k,
\enq
and $\mb_N$ is similarly defined with the $Y$ sample. The vectors  
$\ma_N$ and $\mb_N$ follow multinomial distributions.  
Squared distances between the centers of the cells of $\calc$ are written into the $k\times k$ matrix $D=(D_{j,l})$. For the moment, both the partition $\calc$ and the matrix $D$ are fixed.
In this setting, Bigot et al., \cite{bigot}, prove a Central Limit Theorem for 
$S_\lambda(\ma_N,\mb_N)$, by establishing the Hadamard differentiability of $S_\lambda$ respect to the 
vectors $\ma_N$ and $\mb_N$.  

We want to consider the use of $S_\lambda(\ma_N,\mb_N)$ in a permutation test procedure for the two-sample
problem, adding the flexibility of using a data dependent partition, 
$\calc_N=\{C_{1,N},\dots,C_{k,N}\}$. In particular, we will consider the partition
determined by $k$-means clustering of the joint sample. $k$-means clustering seeks
a collection of centers $\hat\calm=\{\hamu_1,\dots,\hamu_k\}\subset\real^d$ 
that minimize the sum of squared distances from the sample points to the nearest 
vector in $\hat\calm$. Both the sample under consideration and $\real^d$ are 
partitioned into the regions closer to $\hamu_j$ than to any other
$\hamu_l$, for $j,l\leq k$, $l\neq j$. These disjoint regions are convex
polyhedra called Voronoi cells. They are the ``clusters'' forming the partition 
$\calc_N$ produced by the $k$-means algorithm. The center 
of each cluster in the $k$-means algorithm, is the mean of all the data points 
falling in that cluster.  

With $S_\lambda(\ma_N,\mb_N)$ computed as suggested in the previous paragraph, the
calculation of the $p$-value of the statistic by a permutation procedure is performed by the following
steps:  
{\it \begin{itemize}\vskip -0.5 in
\item[(i)] Fix a large integer $B$ (1000 in our simulations) that corresponds to the number of permutations.
Combine the $X$ and $Y$ samples in a joint sample $Z$, of size $n+m$. The Voronoi cells corresponding to the $k$-means clustering applied to the combined
sample are assumed to be known, as well as the value of $S_\lambda(\ma_N,\mb_N)$ for the original samples.
\item[(ii)] For ($b \in 1:B$), do:
\neno (a) Randomly partition $Z$ into samples $X^{(b)}$ and $Y^{(b)}$ of sizes $n$ and $m$, 
respectively. One way of doing this is to reorder $Z$ at random and stating that the resulting first $n$ elements of $Z$ form the sample $X^{(b)}$ and the rest form $Y^{(b)}$.   
\neno (b) Compute the probability vectors $\ma_N^{(b)}$ and $\mb_N^{(b)}$ for the Voronoi cells of the combined sample, with respect to the $X^{(b)}$ and $Y^{(b)}$ samples, respectively.
\neno (c) Compute $S_\lambda(\ma_N^{(b)},\mb_N^{(b)})$. 
\item[(iii)] The approximate permutation $p$-value for the statistic 
$S_\lambda(\ma_N,\mb_N)$ is given by
\[
p\mbox{-value}=\frac{\#\{b\leq B:S_\lambda(\ma_N^{(b)},\mb_N^{(b)})\geq
S_\lambda(\ma_N,\mb_N)\}}{B},
\]
the fraction of times that $S_\lambda$ on the permuted samples is greater or
equal to the original statistic.
\end{itemize}
} 

One relevant computational advantage of using a permutation test in our 
procedure, is that when the data identity is permuted, the clustering structure
associated to $k$-means remains unaltered (since it is associated to the
combined sample), simplifying the calculations required.
This fact would not hold for a classical bootstrap procedure.  

In the two-sample problem, using the partition obtained from $k$-means clustering of 
the joint sample, we can let the matrix $D$ in $S_\lambda(\ma_N,\mb_N)$ 
be the inter center
distance matrix, that is, $D_{j,l}=\|\hamu_j-\hamu_l\|^2$. This will be our {\it basic setting} and the one for which theory is developed in Section 3. Seeking more power for the statistic, we consider the following alternative.    Let $\hamu_{j,X}$ and $\hamu_{j,Y}$ denote, respectively, the 
means of $X$ and $Y$ data falling in cell $C_{j,N}$, for each $j\leq k$. Under the alternative hypotheses, $F\neq G$, we would expect to observe 
the following: (i) For all (or at least some) of the cells, $C_{j,N}$, the 
cell frequencies for the $X$ and $Y$ samples, $a_{j,N}$ and $b_{j,N}$ will differ
noticeably. (ii) For all (or at least some) of the cells, $C_{j,N}$, the $X$ and $Y$
means, $\hamu_{j,X}$ and $\hamu_{j,Y}$ will differ noticeably. To have power
against the general alternative, we want our statistic to be sensitive to both 
conditions. For this purpose, in the definition of $S_\lambda(\ma_N,\mb_N)$ , 
the finite set $\Omega$ is taken as 
\beq\label{omega doble}
\Omega = \{\hamu_{1,X},\dots,\hamu_{k,X},\hamu_{1,Y},\dots,\hamu_{k,Y}\}
\enq 
that is, $L=2k$, and the vectors of frequencies are (re)defined as
\beq\label{freq doble}
\ma_N=(a_{1,N},\dots,a_{k,N},0,\dots,0)'\>\mbox{ and }\>
\mb_N=(0,\dots,0,b_{1,N},\dots,b_{k,N})',
\enq
where $\ma_N$ ends in $k$ zeros and $\mb_N$ begins with $k$ zeros. The matrix 
$D=D_N$ in (\ref{sinkhorn}) is the square distance matrix for the set $\Omega$
given in (\ref{omega doble}). This setting will be called {\it double centers}
in what follows. When one of the samples is absent from one of the
cells, (let us say, for instance that  there is no $Y$ data in $C_{j,N}$), the
corresponding frequency, $b_{j,N}$ is zero, and the center of the absent ($Y$)
sample for the cell, $\hamu_{j,Y}$, is taken as $\hamu_{j',Y}$, where $\hamu_{j',Y}$
is the $Y$ sample center closest to the $X$ center of $C_{j,N}$. 

\section{Theoretical results}
In what follows, a central limit theorem will be established for 
$\hs_\lambda(\ma_N,\mb_N)$ in the {\it basic setting} defined towards the end of Section 2. The proof, relatively elementary, uses the fact that the number of iterations required by the Sinkhorn algorithm remains bounded when computed on probability vectors on a fixed number of cells. 

Let $\mu^*=(\mu_1^*,\dots,\mu_k^*)\in\real^{kd}$ be the vector of centers corresponding to the population solution of the $k$ means functional for the distribution $P$ of the $X$ sample. $\mu^*$ is assumed to be unique.
Denote by $\cals$ the support of $P$. We assume
that $\cals$ is a bounded set in $\real^d$. 
Let $\calk$ be a compact neighborhood of $\mu^*$, of small diameter. To each 
$\mu=(\mu_1,\dots,\mu_k)\in \calk$, assign the partition 
$\calc_\mu=\calc=\{C_1(\mu),\dots,C_k(\mu)\}$ of Voronoi cells associated to the centers in
$\mu$. Consider the canonical vectors $e_1,\dots,e_{k-1}$ in $\real^{k-1}$: $e_l$ has 
all coordinates equal to zero, except for a 1 in the $l$-th coordinate.
 Define the function \linebreak  $h_\mu:\real^d\rightarrow\real^{k-1}$ given by \[ h_\mu(x)=e_l \mbox{ if } x\in C_l(\mu),\>\>\mbox{ for }l\leq k-1\>\>\mbox{ and } h_\mu(x)=\mathbf{0}\mbox{ if } x\in C_k(\mu).\] 
The $h_\mu(X_i)$ are the multivariate ``Bernoulli'' vectors whose average produce
the first $k-1$ coordinates of the probability vector $\ma_N$ that corresponds to 
the partition associated to $\mu$. 

First, we have the following Lemma. 

\begin{lem}\label{lema1}
To the setting of the previous paragraph, add the assumption that the probability distribution $P$ admits a bounded, continuous density. Then, the class of functions 
$\calh=\{h_\mu:\mu\in\calk\}$ is a $P$-Donsker class, with covering $L^2$ number,
$N_2(\eps,\calh,P)$, satisfying
\[
N_2(\eps,\calh,P)\leq A\left(\frac{1}{\eps}\right)^{2kd},
\]
for a positive constant $A$. 
\end{lem}
\proof{For $h_\mu(X)- h_{\mu'}(X)\neq \mathbf{0}$ to hold, there must exist an index $r$,
such that $X$ belongs to the symmetric difference 
\[C_r(\mu)\Delta C_r(\mu')=(C_r(\mu)\setminus C_r(\mu'))\cup
(C_r(\mu')\setminus C_r(\mu)).\] Thus,
$\rho^2_P(h_\mu,h_{\mu'}):=P\left(\|h_\mu-h_{\mu'}\|^2\right)\leq 2 \sum_{r}P(C_r(\mu)\Delta C_r(\mu')),$ (since \linebreak $\|h_\mu-h_{\mu'}\|^2\leq\|e_r\|^2+\|e_{r'}\|^2=2$). To bound the
covering number of $\calh$ with respect to $\rho_P$, we will bound the 
change in $P(C_r(\mu')\setminus C_r(\mu))$ in terms of the norm $\|\mu-\mu'\|$,
when $\mu'$ is obtained by a small change in $\mu$. This
reduces, in turn, to considering the effect on one face of the Voronoi cell $C_r$.
Without loss of generality, assume that $r=1$ and we are considering the face
between cells $C_1$ and $C_2$. Again, without loss of generality, assume that
$\mu_1$ and $\mu_2$ are placed along the $x_1$ axis and $\mathbf{0}$ is the middle
point between them. That is, for some $\Delta>0$,
\[
\mu_1=(\Delta,0,\dots,0)\>\mbox{ and }\> \mu_2=(-\Delta,0,\dots,0).
\]
Then, in the face between cells $C_1$ and $C_2$, $C_1$ is chosen if the first
coordinate is positive ($x_1>0$). Suppose $\mu_2$ changes to
$\mu'_2=(-\Delta+\delta_1,\delta_2,\dots,\delta_d)$, while $\mu_1$ does
not change. The midpoint between $\mu_1$ and $\mu'_2$ is now
$(1/2)(\delta_1,\delta_2,\dots,\delta_d)$ and the half-space for choosing $C'_1$ 
over $C'_2$ becomes 
\[
(2\Delta-\delta_1)x_1-\delta_2x_2\dots -\delta_dx_d>
(2\Delta-\delta_1)\frac{\delta_1}{2}-\frac{\delta_2^2}{2}\dots 
-\frac{\delta_d^2}{2},
\]
for $x=(x_1,\dots,x_d)\in\real^d$. With respect to the $1$-$2$ face, the 
condition \linebreak $x\in C'_1\setminus C_1$ reduces to
\beq\label{cond x_1}
0\geq 2\Delta x_1 > \delta_1x_1+\dots +\delta_dx_d
+\Delta{\delta_1}-\frac{\delta_1^2}{2}-\frac{\delta_2^2}{2}\dots 
-\frac{\delta_d^2}{2},
\enq
which, by the bounded support condition, can be written as \linebreak
$0\geq x_1> \mbox{O}(\|\mu_2-\mu_2'\|)$. Using the continuity of the density, it follows that the probability of condition (\ref{cond x_1}) is 
$\mbox{O}(\|\mu_2-\mu_2'\|)$. A similar analysis for every face of the cell 
$C_1$, leads to
\[
\Pr(C'_1\setminus C_1)=\mbox{O}(\|\mu-\mu'\|),
\]
a bound that holds for every Voronoi cell and we get
\beq
\rho_P(h_\mu,h_{\mu'})\leq \mbox{O}(\|\mu-\mu'\|^{1/2}).
\enq 
Since the compact set 
$\calk$ has a finite diameter in $\real^{kd}$, it can be covered, with respect to
Euclidean distance, within distance $\epsilon^2$, by an array of points,
$\calz_\eps$, of 
cardinality $\mbox{O}((1/\epsilon)^{2kd})$: For each $\mu\in\calk$,
there is a $\mu'\in \calz_\eps$, such that
\beq\label{cota1}
P\left(\|h_\mu-h_{\mu'}\|^2\right)\leq \eps^2.
\enq
The rest of the proof follows an argument similar to that of Example 19.7 in \cite{van der vaart}, using
covering number instead of bracketing covering number. 
\qed}

The computational complexity of computing the Sinkhorn discrepancy between two probability vectors of length $L$, is $\mbox{O}(L^2)$ (see \cite{awr17}, \cite{cuturi}). In our context, when the number of clusters in the $k$-means procedure is fixed, this means that the computational complexity is $\mbox{O}(1)$.
In particular, the number of iterations that the Sinkhorn procedure performs,
before reaching a stopping criteria, is bounded. In
our experiments, we have observed that the distribution over 1,000 permutation
replicas, of the number of Sinkhorn iterations 
depends, to some extent,  on the (mixture) distribution of the joint sample, and 
is typically supported on three or four consecutive integers, 7, 8, 9 and 10,
for example. It is not difficult to determine, experimentally, a very likely
upper bound for the number of iterations in each case. 

For the probability vectors, $\ma,\mb$ and the cost matrix $D$ appearing in 
(\ref{sinkhorn}), let $K=\exp(-D/\lambda)$, 
with operations performed elementwise. The Sinkhorn iteration, $R:\real^L\times\real^L\rightarrow\real^L\times\real^L$,
is given by the transformation

\beq\label{S iter}
R(u,v;\ma,\mb)=\left(\frac{\ma}{Kv},\frac{\mb}{K^tu}\right),
\enq
where, again, operations are performed elementwise. Formula (\ref{S iter}) is 
used to compute new values of $u$ and $v$, from the previous ones, until convergence
according to a ``negligible change'' criteria. The values obtained define the
solution of (\ref{sinkhorn}) since $T^\lambda=\mbox{diag}(u)K\mbox{diag}(v)$ achieves
the minimum of (\ref{sinkhorn}) (approximately, due to the stopping criteria).
In our basic setting (as defined towards the end of Section 2), using vectors
of positive coordinates as initial values for $u$ and $v$, the coordinates of
these vectors remain positive throughout the iterations and all the entries of 
the resulting $T^\lambda$ are positive too (see \cite{cuturi2014}). By composition of differentiable maps, it follows that the transformations leading to the final
values of $u$ and $v$, to $T^\lambda$ and to the optimum value 
$S_\lambda(\ma,\mb)=\left\langle D,T^\lambda\right\rangle
-\lambda H(T^\lambda)$ are differentiable with respect to the vectors $\ma$, $\mb$ 
and to the vector of centers, $\mu=(\mu_1,\dots,\mu_k)$, with partial derivatives 
uniformly bounded on the set $\calk$ of the previous Lemma.

Let us introduce some additional elements required for our Central Limit Theorem.
Let $\ma^*=(a^*_1,\dots,a^*_k)$ be the vector of
cell probabilities associated to the vector of population centers, 
$\mu^*=(\mu^*_1,\dots,\mu^*_k)$. Then, the covariance matrix for the Bernoulli vector $h_{\mu^*}(X_i)$ is
\beq\label{cova}
\Sigma=\Sigma(\mu^*)=\left(\begin{array}{cccc}
a^*_1(1-a^*_1) & -a^*_1a^*_2 & \cdots & -a^*_1a^*_{k-1} \\
-a^*_2a^*_1 & a^*_2(1-a^*_2) & \cdots & -a^*_2a^*_{k-1} \\
&\cdots &\cdots & \\
-a^*_{k-1}a^*_1 & -a^*_{k-1}a^*_2 & \cdots & a^*_{k-1}(1-a^*_{k-1})
\end{array}
\right)
\enq
Let $I$ be an upper bound for the number of iterations required by the Sinkhorn procedure, for the distribution generating the two sample data under the
null hypothesis and for the value of $k$ used in the $k$ means clustering procedure.  
Assume Sinkhorn divergences are computed using always $I$ iterations.
Assume also that the sample sizes, $n$ and $m$, satisfy $\>\> n/N\to\alpha\in(0,1)$ 
as $n,m\to\infty$. 

Let ${\ma}_N$ and ${\mb}_N$ denote the sample frequency vectors, for the $X$ and $Y$ samples, for the partition associated to $\mu^*$.
Let $\hat{\ma}_N$ and $\hat{\mb}_N$ denote the  
sample frequencies associated to the partition produced by the $k$ means 
clustering (as described in the lines following (\ref{freq})) and let $\ma_N^*$ be the population probability vector for the same data dependent cells. To remove degenerate distributions, we
will also need the truncated frequency vectors $\hat\ma_{N,t}$ and
$\ma^*_t$ defined as the
vectors of the first $k-1$ coordinates of $\hat\ma_{N}$ and $\ma^*$, 
respectively. Define, similarly, $\hat\mb_{N,t}$ and
$\mb^*_t$.

Let ${\cal H}_I(\ma,\mb)$ be the transformation defined by the following steps:
(i) For positive initial vectors, $u_0$ and $v_0$, apply transformation 
(\ref{S iter}) $I$ times, to reach the final value of the pair $(u,v)$, (ii) with 
these $u$ and $v$, compute \linebreak $T^\lambda=\mbox{diag}(u)K\mbox{diag}(v)$, and (iii) 
apply formula (\ref{sinkhorn_3}) to obtain the value of $\hs_\lambda(\ma,\mb)$.
By the discussion after (\ref{S iter}), ${\cal H}_I$ is a differentiable function
of  $\ma$ and $\mb$. Denote by ${J}_I$ the gradient of ${\cal H}_I$ with respect
to the first $k-1$ coordinates of $\ma$ and $\mb$, written as a row vector. In this setting, we have
\begin{theo}\label{clt1}
 Let $N_{\mbox{e}}=nm/N$. 
Under the null hypothesis and the assumptions listed for Lemma \ref{lema1},
\beq\label{limit1}
\sqrt{N_{\mbox{e}}}(\hs_\lambda(\hat{\ma}_N,\hat{\mb}_N)-\hs_\lambda({\ma}^*_N,{\ma}^*_N)) \cid N(0,J_I C J_I')
\enq 
where $C$ is given by
\[
C=\left(
\begin{array}{cc}
(1-\alpha)\Sigma & \mathbf{0} \\
\mathbf{0} & \alpha\Sigma 
\end{array}
\right)
\]
as $n,m\to\infty$, and the $\mathbf{0}$s in the definition of $C$ are blocks of
the appropriate size.
\end{theo}

\proof{For the limiting centers in $\mu^*$ and the corresponding cells, 
by the Central Limit Theorem for the multinomial distribution, under the null hypothesis we have 
\beqa\label{clt1a}
\sqrt{n}(\ma_{N,t}-\ma^*_{t})\cid N(\mathbf{0}, \Sigma) & \mbox{ and } \nonumber \\
\sqrt{m}(\mb_{N,t}-\ma^*_{t})\cid N(\mathbf{0}, \Sigma) & 
\enqa
Multiplying the lines of (\ref{clt1a}) by $\sqrt{m/N}$ and $\sqrt{n/N}$ respectively, and using the independence of the two samples, we get 
\beq\label{clt1b}
\nu_N(\mu^*):=\sqrt{N_e}((\ma_{N,t},\mb_{N,t})-(\ma^*_{t},\ma^*_{t}))\cid N(\mathbf{0}, C).
\enq
By the main theorem in Section 3, of 
\cite{pollard82}, $\hamu\cip \mu^*$, as $n,m\to\infty$. Let $\nu_N(\hat\mu)=\sqrt{N_e}((\hat\ma_{N,t},\hat\mb_{N,t})-(\ma^*_{N,t},\ma^*_{N,t}))$. By Dudley's Asymptotic Equicontinuity Condition, implied by Lemma 1, we have $\nu_N(\hat\mu)-\nu_N(\mu^*)=\och(1)$. Therefore,
\beq\label{clt1c}
\nu_N(\hat\mu):=\sqrt{N_e}((\hat\ma_{N,t},\hat\mb_{N,t})-(\ma^*_{N,t},\ma^*_{N,t})) \cid N(\mathbf{0}, C).
\enq
Apply the function $\mathcal{H}_I$ to both terms in (\ref{clt1c}). Notice that $ \mathcal{H}_I(\hat{\ma}_{N,t},\hat{\mb}_{N,t})$  and ${\mathcal H}_I(\ma^*_{N,t},\ma^*_{N,t})$ are equal to $\hs_\lambda(\hat{\ma}_N,\hat{\mb}_N)$ and $\hs_\lambda({\ma}^*_N,{\ma} 
^*_N)$ respectively. By the usual multivariate version of the Delta Method, 
\beq
\sqrt{N_{\mbox{e}}}(\hs_\lambda(\hat{\ma}_N,\hat{\mb}_N)-\hs_\lambda({\ma}^*_N,{\ma}
^*_N)) \cid N(0,J_I C J_I').
\enq 
\qed}

Similar results hold for $S_\lambda(\hat{\ma}_N,\hat{\mb}_N)$ and
$\os_\lambda(\hat{\ma}_N,\hat{\mb}_N)$, and can be proved with slight 
redefinitions of the transformation $\mathcal{H}_I$ in each case. It is also
possible to obtain a similar result for the distribution of 
$\hs_\lambda(\hat{\ma}_N,\hat{\mb}_N)$ under the alternative hypothesis.

The gradient vector $J_I$ in Theorem \ref{clt1} can be consistently estimated
by computing $\hs_\lambda$ on small perturbations of the coordinates of 
$\hat{\ma}_{N,t}$ and $\hat{\mb}_{N,t}$. Since the matrix $C$ is also estimable
from the sample (and the $k$-means clustering cells), 
the limiting distribution in (\ref{limit1}) can be 
specified in an approximate (consistent) way. On this regard, see the experiment
described at the end of next Section.

\section{Performance evaluation on simulated data}

In this section, we will compare $S_\lambda(\ma,\mb)$, $\os_\lambda(\ma,\mb)$ and $\hs_\lambda(\ma,\mb)$, in terms of statistical power and computational cost, among them, and also against Wasserstein's statitic, $W(\ma,\mb)$ and Schilling's statistic (\cite{schilling}), known for being one of the most powerful tests for the non-parametric two sample problem. 
$S_\lambda(\ma,\mb)$, $\os_\lambda(\ma,\mb)$, $\hs_\lambda(\ma,\mb)$ are all implemented with {\it double centers} in the partition
cells (as described in the previous section). All the statistics were implemented
in the R Language by the authors as permutation tests, using functions from package {\it Barycenter} for the Sinkhorn statistics.

\subsection{Simulated data description}

In order to compare the statistics previously mentioned, we consider three different continuous distribution scenarios for the data:
\begin{enumerate}
	\item \textbf{Gaussian Multivariate Distributions:} In this case, we use \linebreak $5$-dimensional Gaussian distributions. $X$ is generated from a Standard Gaussian Distribution and $Y$ is generated from the same distribution of $X$ (null hypothesis case) or from one of the following alternative distributions:
	\begin{itemize}
		\item A Gaussian distribution with mean \begin{small} $\mu_1=(0.05, 0.01, -0.05, 0, 0.101)$ \end{small} and covariance matrix $I_5$.
		\item A Gaussian distribution with mean \begin{small} $\mu_2=(0.11, 0.022, -0.011, 0, 0.222)$ \end{small}  and covariance matrix $I_5$.
		\item A Gaussian distribution with mean \begin{small} $\mu_3=(0.5, 0.1, -0.5, 0, 1.01)$ \end{small} and covariance matrix $I_5$.
		\item A Gaussian distribution with mean \begin{small} $\mu=(0, 0, 0, 0, 0)$ \end{small} and covariance matrix 
		\begin{small}
			\[
			\Sigma_1=\left[\begin{array}{ccccc}
				1.065  & 0.044  & -0.036 & 0.01   & 0.019  \\
				0.044  & 1.081  & 0.006  & 0.023  & -0.016 \\
				-0.036 & 0.006  & 1.066  & -0.016 & -0.024 \\
				0.01   & 0.023  & -0.016 & 1.046  & -0.026 \\
				0.019  & -0.016 & -0.024 & -0.026 & 1.039
			\end{array}\right].
			\] 
		\end{small}
		\item A Gaussian distribution with mean \begin{small} $\mu=(0, 0, 0, 0, 0)$ \end{small} and covariance matrix 
		\begin{small}
			\[
			\Sigma_2=\left[\begin{array}{ccccc}
				1.1475 &  0.066  & -0.054 &   0.015  &  0.0285  \\
				0.066  &  1.1715 &  0.009 &   0.0345 & -0.024   \\
				-0.054  &  0.009  &  1.149 &  -0.024  & -0.036   \\
				0.0150 &  0.0345 & -0.024 &   1.1190 & -0.039   \\
				0.0285 & -0.024  & -0.036 &  -0.039  &  1.1085
			\end{array}\right].
			\]
		\end{small}
		\item A Gaussian distribution with mean \begin{small} $\mu=(0, 0, 0, 0, 0)$ \end{small} and covariance matrix 
		\begin{small}
			\[
			\Sigma_3=\left[\begin{array}{ccccc}
				1.65  & 0.44  & -0.36 & 0.1   & 0.19  \\
				0.44  & 1.81  & 0.6   & 0.23  & -0.16 \\
				-0.36 & 0.6   & 1.66  & -0.16 & -0.24 \\
				0.1   & 0.23  & -0.16 & 1.46  & -0.26 \\
				0.19  & -0.16 & -0.24 & -0.26 & 1.39
			\end{array}\right].
			\]
		\end{small}		
	\end{itemize}
	The choice of means and covariance matrices in the alternatives just presented
		seeks to achieve a gradual departure from the null hypothesis (the standard
		Gaussian distribution) in mean or in covariance matrix.
	\item \textbf{Bounded Burr Distributions:} Given $E_1,\, \ldots ,\, E_d$ i.i.d. variables with the Exponential($\theta=1$) distribution, and $G$ a
	random variable independent of the $E_i$, with 
	$G\sim\mbox{Gamma}\left(\alpha,\beta=1\right)$,
	for $\alpha>0$.
	Let  
	\[
	B_i = \dfrac{E_i}{G}, \quad 1\leq i \leq d.
	\]
	and define the random vector $B=(B_1,\dots,B_d)$. The distribution of $B$ is what we call the Bounded Burr. It appears in \cite{johnson}, as an auxiliary family. One of the characteristics of this family, which motivates its 
		use in our experiments, is that it has compact support contained in
		the $d$-dimensional hypercube where, as the parameter $\alpha$ increases above 1, the data become more uniformly distributed on the 
		hypercube, while, as $\alpha$ decreases towards zero, the support of the distribution  shrinks towards a region around a diagonal of the hypercube. This
		type of support offers a clear variation from the Gaussian distribution.
	
	In this setting, $X$ is generated from a Bounded Burr distribution with parameter $\alpha=1$ and, for the alternatives, $Y$ is generated from a Bounded Burr distribution with parameters $\alpha_1=1$, $\alpha_2=1.25$, $\alpha_3=1.5$, \linebreak $\alpha_4=1.75$,  $\alpha_5=2$, $\alpha_6=0.75$, $\alpha_7=0.5$.  
	
	The fact that this distribution depends only on the parameter $\alpha$, makes it  simpler to generate a number of alternatives controlling the grade of difference with respect to the null case.
	
	\item \textbf{Dirichlet Distribution:} This family of continuous distributions,
	described in \cite{kbj}, is a multivariate generalization of the Beta distribution. The support of the Dirichlet data is restricted to the simplex of $d$-dimensional probability vectors, thus offering a clearly different
		situation from the distributions in our previous settings. The distribution is defined by a concentration parameter $\beta\in\mathbb{R}^d$. On this scenario, $X$ (and $Y$ under the null hypothesis) will be generated from a Dirichlet distribution with parameter 
	$\beta_1=(1,1,1,1,1)$, and, for the alternatives $Y$ will be generated from the Dirichlet distribution with parameters:
	\begin{itemize} 
		\item $\beta_2=c(1.5,1.5,1.5,1.5,1.5)$. 
		\item $\beta_3=(0.5,0.5,0.5,0.5,0.5)$.
		\item $\beta_4=(1.015, 0.95, 1.043, 0.975, 0.98)$. 
		\item $\beta_5=(1.27, 0.59, 1.55, 1.23, 0.36)$. 
		\item $\beta_6=(2,2,2,2,2)$.
	\end{itemize} 
	
\end{enumerate} 

In our first experiments, we use samples of sizes $n=m=500$ in all cases. All the
statistics are used as permutation statistics. For this purpose, 1000 random permutations of the combined sample are used. For the Schilling test, $K=4$ neighbors are used. Figures 1 to 4, show histograms of the p-values obtained for $S_\lambda(\ma,\mb)$ and Schilling's statistic in the Gaussian distribution case for the null distribution and alternatives $\mu=\mu_1$, $\mu=\mu_2$ and $\Sigma=\Sigma_2$. 

\begin{figure}[h!]
	\centering
	\includegraphics[scale=0.313]{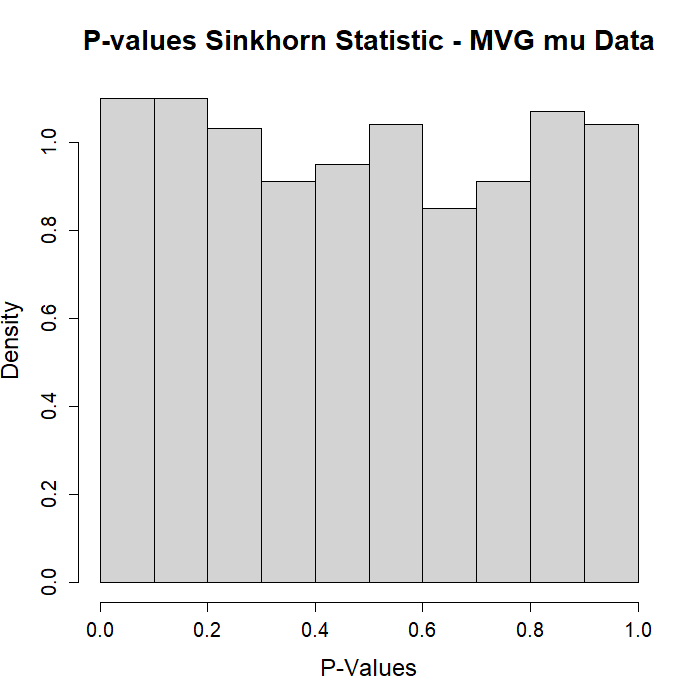}
	\includegraphics[scale=0.313]{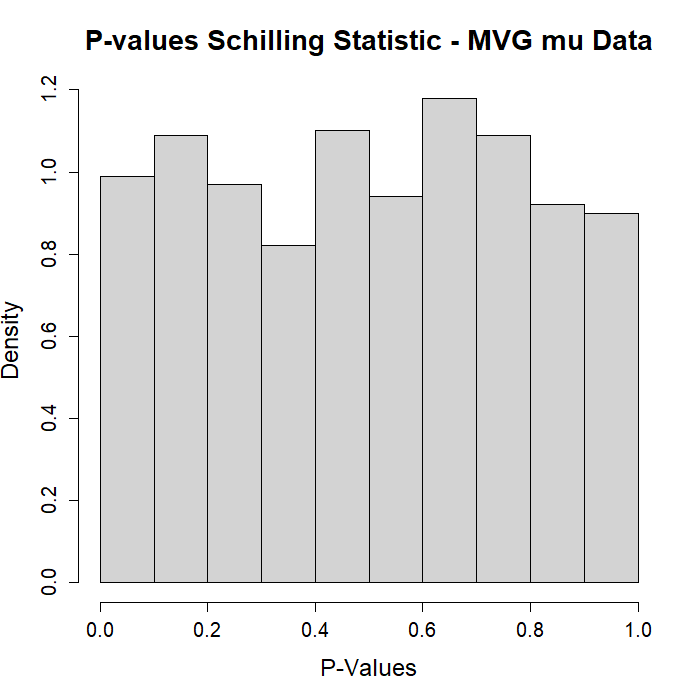}
	\caption{Histograms of p-values computed with statistic $S_\lambda(a,b)$ (Left)
		and Schilling (Right) in the case of MVG data under the null hypothesis.}
	\label{fig:pval0}
\end{figure}

\begin{figure}[h!]
	\centering
	\includegraphics[scale=0.313]{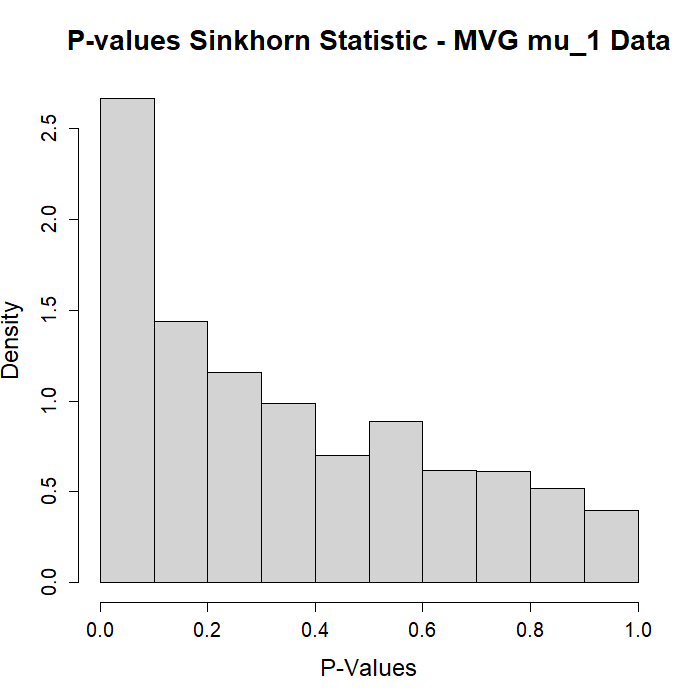}
	\includegraphics[scale=0.313]{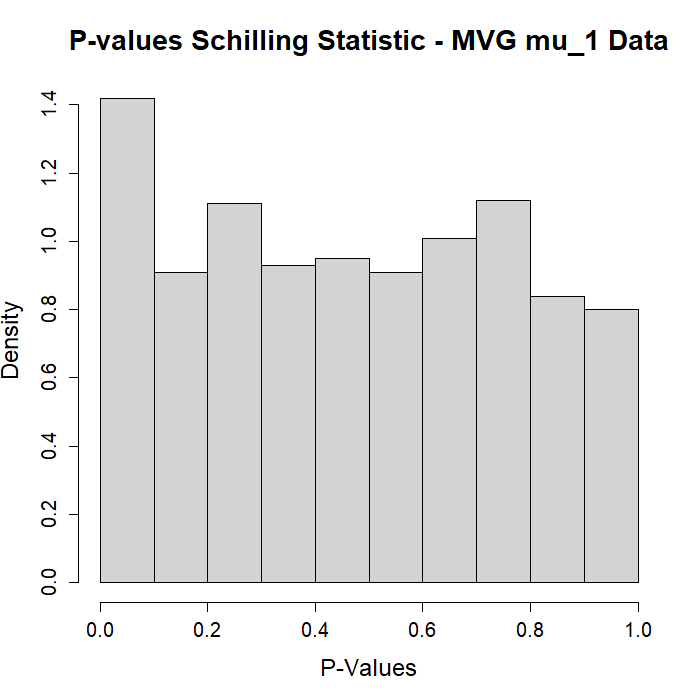}
	\caption{Histograms of p-values computed with statistic $S_\lambda(a,b)$ (Left)
		and Schilling (Right) in the case of MVG data under the alternative hypothesis where $Y$ is distributed with mean $\mu_1$.}
	\label{fig:pval1}
\end{figure}

\begin{figure}[h!]
	\centering
	\includegraphics[scale=0.313]{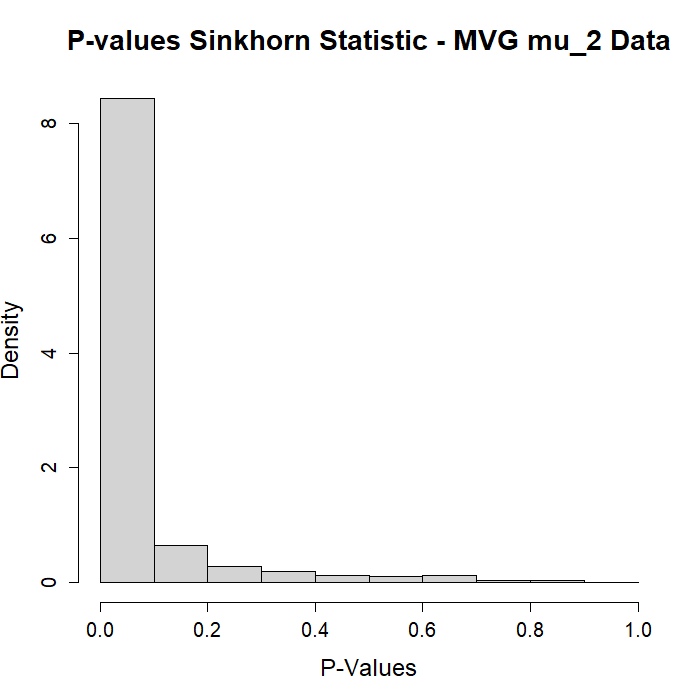}
	\includegraphics[scale=0.313]{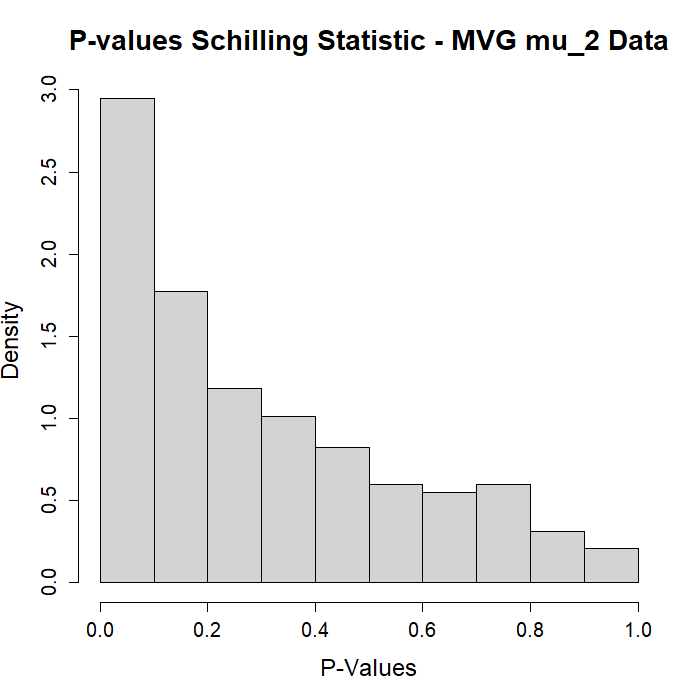}
	\caption{Histograms of p-values computed with statistic $S_\lambda(a,b)$ (Left)
		and Schilling (Right) in the case of MVG data under the alternative hypothesis where $Y$ is distributed with mean $\mu_2$.}
	\label{fig:pval2}
\end{figure}

As should be expected, in the null hypothesis case, the distribution of the $p$-values resembles the uniform distribution. For both of the mean alternatives in Figures 2 and 3, the permutation test based on $S_\lambda(\ma,\mb)$ displays better power that Schilling's test, with a $p$-value distribution more concentrated near zero. The same behavior is observed in the covariance matrix alternative (Figure 4).

\begin{figure}[h!]
	\centering
	\includegraphics[scale=0.313]{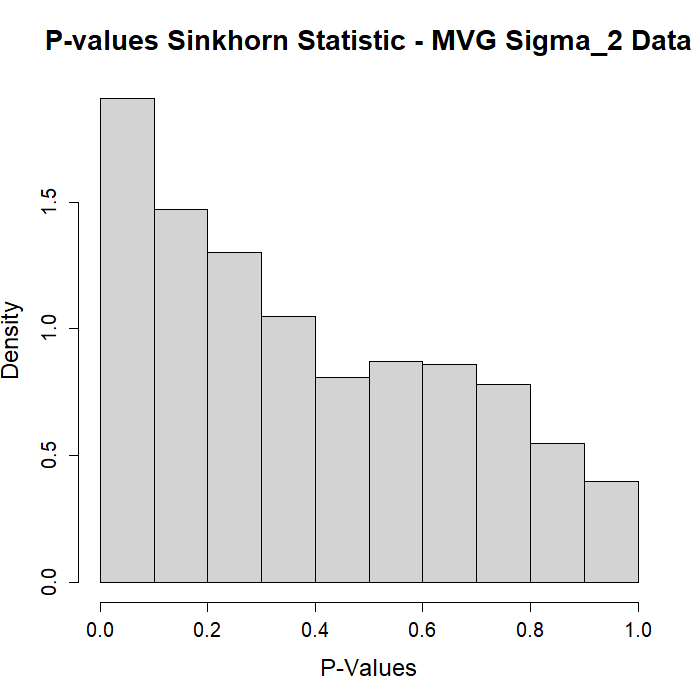}
	\includegraphics[scale=0.313]{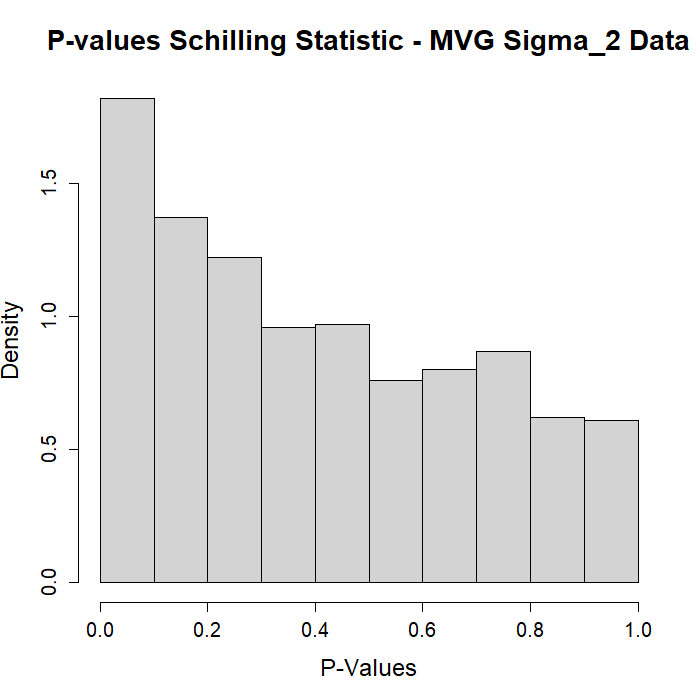}
	\caption{Histograms of p-values computed with statistic $S_\lambda(a,b)$ (Left)
		and Schilling (Right) in the case of MVG data under the alternative hypothesis where $Y$ is distributed with covariance matrix $\Sigma_2$.}
	\label{fig:pval5}
\end{figure}

\subsection{Power comparison}

Table 1 contains a more detailed power comparison, at the 5\% level of the permutation tests based on the different statistics considered.

\begin{table}[]
	\centering
	\resizebox{0.9\textwidth}{!}{%
	\begin{tabular}{|c|c|ccccc|}
		\hline
		\multirow{2}{*}{ } & \multirow{2}{*}{ Data Set } & \multirow{2}{*}{ $ Schil. $ } & \multirow{2}{*}{ $ W(\ma,\mb) $ }
		& \multirow{2}{*}{ $\, S_\lambda(\ma,\mb) \,$ } & \multirow{2}{*}{$\, \hs_\lambda(\ma,\mb) \,$}
		& \multirow{2}{*}{$\, \os_\lambda(\ma,\mb) \,$ } \\[-5pt]
		&  &  &  &  &  & \\
		\hline
		\multirow{7}{*}{\rotatebox[origin=c]{90}{MVG}}  
		& $\mu$      & 5.8  & 5.5  & 4.0  & 5.2  & 5.0  \\
		& $\mu_1$    & 6.7  & 12.6 & 16.5 & 14.9 & 15.0 \\
		& $\mu_2$    & 20.4 & 51.2 & 73.3 & 70.3 & 70.7 \\
		& $\mu_3$    & 100  & 100  & 100  & 100  & 100  \\
		& $\Sigma_1$ & 6.9  & 3.5  & 7.7  & 6.3  & 6.5  \\
		& $\Sigma_2$ & 10.8 & 2.5  & 11.9 & 10.0 & 9.2  \\
		& $\Sigma_3$ & 99.7 & 22.7 & 100  & 99.9 & 99.7 \\ 
		\cline{2-7}
		&$avg.\,time$ &1.42  & 1.57 & 0.71 & 0.69 & 0.84 \\
		\hline
		\multirow{7}{*}{\rotatebox[origin=c]{90}{Bounded Burr}} 
		& $\alpha_1$ & 4.9  & 4.3  & 5.5  & 5.3  & 4.9  \\
		& $\alpha_2$ & 5.3  & 7.8  & 11.7 & 10.0 & 13.0 \\
		& $\alpha_3$ & 8.5  & 15.0 & 30.9 & 24.1 & 33.6 \\
		& $\alpha_4$ & 13.4 & 32.5 & 57.9 & 53.5 & 61.9 \\
		& $\alpha_5$ & 19.1 & 50.9 & 78.0 & 78.3 & 81.8 \\
		& $\alpha_6$ & 9.3  & 11.2 & 12.2 & 12.2 & 14.2 \\
		& $\alpha_7$ & 50.9 & 71.2 & 87.1 & 74.4 & 80.1 \\ 
		\cline{2-7}
		&$avg.\,time$ & 1.22 & 1.56 & 0.79 & 0.63 & 0.75 \\
		\hline
		\multirow{7}{*}{\rotatebox[origin=c]{90}{Dirichlet}} 
		& $\beta_1$ & 5.5  & 6.3  & 4.5  & 5.7  & 4.5  \\
		& $\beta_2$ & 83.5 & 62.7 & 83.9 & 91.4 & 33.9 \\
		& $\beta_3$ & 100  & 13.0 & 100  & 100  & 100  \\
		& $\beta_4$ & 7.7  & 7.8  & 9.5  & 12.7 & 14.4 \\
		& $\beta_5$ & 100  & 100  & 100  & 100  & 100  \\
		& $\beta_6$ & 100  & 99.9 & 100  & 100  & 98.7 \\
		\cline{2-7}
		&$avg.\,time$ & 1.38 & 1.60 & 0.89 & 0.68 & 0.77 \\
		\hline 
	\end{tabular}}
	\caption{Empirical power of the tests, measured in \% for simulated experiments on $k=10$ clusters and $\lambda=1$, and $K=4$ nearest neighbors for the Schilling test. \textit{avg. time} is measured in seconds.}
	\label{fig: Table1}
\end{table}

In  Table \ref{fig: Table1}, (in the first row corresponding to each distribution)
we see that all the permutation tests respect adequately the type I error level 
under the null hypothesis, in the different distributional scenarios. In terms of 
power, the numbers displayed in this table can be summarized as follows: In the case 
of the MVG distributions, for mean alternatives, the three tests based on the
Sinkhorn divergence perform better that Schilling's statistic and the test based
on the Wasserstein distance, while in the case of MVG distributions, for covariance
alternatives, Schilling's statistic and the Sinkhorn tests perform similarly and
slightly better than $W(\ma,\mb)$.
In the case of data with the Bounded Burr distribution, it is clear again that the
three Sinkhorn tests behave similarly and display more power that the test of Schilling and the one based on the Wasserstein distance. Finally, for data with
the Dirichlet distribution, $\hs_\lambda(\ma,\mb)$ exhibits the best power, 
closely followed by $S_\lambda(\ma,\mb)$ and Schilling's statistic, while the other
two statistics fall behind against some of the alternatives. 

As for computational cost, in Table 1, the ``avg. time'' rows show the average
time, in seconds, required to compute a $p$-value (involving 1,000 permutations
and statistic calculations) in each of the scenarios considered. We see that,
in all cases, the calculations required for the Sinkhorn statistics, as
permutation tests, are less than a second, while the statistic based
on $W(\ma,\mb)$ requires about twice as long, and the Schilling statistic
requires an amount of time in between those for the Sinkhorn and
Wasserstein statistics. For the sample sizes, dimension and number of
clusters in the $k$-means procedure, the computational cost of the permutation
tests based on the Sinkhorn divergence are quite manageable.

The following subsections evaluate variations in performance of the Sink\-horn's divergence based statistics, as we consider changes in the value of the parameter
$\lambda$, number of clusters in the $k$-means algorithm, and the effect of changes
in the dimension of the data. 

\subsubsection{Variations on the $\lambda$ parameter.}

In the context of the MVG distributed data of Table \ref{fig: Table1}, we let $\lambda$ take values $0.5$, $0.75$, $1$, $5$, $10$ and $50$, making the optimization problem move from near the Wasserstein distance to near a ``pure entropy''  problem. Table \ref{fig: Table2} shows the powers estimated in this experiment. As $\lambda$ increases from $0.5$ up to $\lambda=5.0$, the power for all versions of the statistic increase noticeably against both the mean and the covariance alternatives. When $\lambda$ goes beyond $5.0$ power figures against the covariance alternatives begin to decrease, and the improvement against the mean alternatives becomes less significant, suggesting that, at least for Gaussian data, $\lambda=5.0$ is a good
choice for the parameter.

\begin{table}[]
	\centering
	\resizebox{0.9\textwidth}{!}{%
	\begin{tabular}{|l|ccc|ccc|}
		\hline
		& \multicolumn{3}{c|}{$\lambda=0.5$} & \multicolumn{3}{c|}{$\lambda=0.75$}  \\ \cline{2-7} 
		& \multirow{2}{*}{ $\, S_\lambda(\ma,\mb) \,$ }  & \multirow{2}{*}{$\, \hs_\lambda(\ma,\mb) \,$}
		& \multirow{2}{*}{$\, \os_\lambda(\ma,\mb) \,$ } & \multirow{2}{*}{ $\, S_\lambda(\ma,\mb) \,$ } & \multirow{2}{*}{$\, \hs_\lambda(\ma,\mb) \,$}    & \multirow{2}{*}{$\, \os_\lambda(\ma,\mb) \,$ } \\[-5pt]
		&  &  &  &  &  &  \\
		\hline
		$\mu$     & 4.9  & 5.8  & 4.8  & 5.3  & 5.3  & 4.4  \\
		$\mu_1$   & 9.8  & 16.1 & 18.3 & 17.7 & 17.0 & 17.4 \\
		$\mu_2$   & 50.7 & 68.2 & 66.7 & 79.1 & 70.2 & 68.1 \\
		$\mu_3$   & 100  & 100  & 100  & 100  & 100  & 100  \\
		$\Sigma_1$& 6.2  & 5.8  & 5.3  & 7.3  & 6.7  & 5.8  \\
		$\Sigma_2$& 12.1 & 9.0  & 9.0  & 12.4 & 8.4  & 7.6  \\
		$\Sigma_3$& 99.8 & 99.8 & 99.9 & 100  & 99.9 & 99.7 \\ 
		\hline
		$avg.\,time$ & 0.96 & 0.80 & 2.30 & 0.84 & 0.70 & 1.06 \\
		\hline\hline
		& \multicolumn{3}{c|}{$\lambda=1$} & \multicolumn{3}{c|}{$\lambda=5$}  \\ \cline{2-7} 
		& \multirow{2}{*}{ $\, S_\lambda(\ma,\mb) \,$ }  & \multirow{2}{*}{$\, \hs_\lambda(\ma,\mb) \,$}
		& \multirow{2}{*}{$\, \os_\lambda(\ma,\mb) \,$ } & \multirow{2}{*}{ $\, S_\lambda(\ma,\mb) \,$ } & \multirow{2}{*}{$\, \hs_\lambda(\ma,\mb) \,$}    & \multirow{2}{*}{$\, \os_\lambda(\ma,\mb) \,$ } \\[-5pt]
		&  &  &  &  &  &  \\ 
		\hline
		$\mu$     & 4.0  & 5.2  & 5.0  & 4.1  & 5.2  & 4.6  \\
		$\mu_1$   & 16.5 & 14.9 & 15.0 & 19.0 & 20.6 & 20.9 \\
		$\mu_2$   & 73.3 & 70.3 & 70.7 & 86.1 & 85.4 & 84.3 \\
		$\mu_3$   & 100  & 100  & 100  & 100  & 100  & 100  \\
		$\Sigma_1$& 7.3  & 6.3  & 6.5  & 7.7  & 7.0  & 7.4  \\
		$\Sigma_2$& 11.9 & 10.0 & 9.2  & 16.3 & 16.1 & 19.7 \\
		$\Sigma_3$& 100  & 99.9 & 99.7 & 100  & 100  & 100  \\ 
		\hline
		$avg.\,time$ & 0.71 & 0.69 & 0.84 & 0.77 & 0.64 & 0.73 \\
		\hline\hline
		& \multicolumn{3}{c|}{$\lambda=10$} & \multicolumn{3}{c|}{$\lambda=50$}  \\ \cline{2-7} 
		& \multirow{2}{*}{ $\, S_\lambda(\ma,\mb) \,$ }  & \multirow{2}{*}{$\, \hs_\lambda(\ma,\mb) \,$}
		& \multirow{2}{*}{$\, \os_\lambda(\ma,\mb) \,$ } & \multirow{2}{*}{ $\, S_\lambda(\ma,\mb) \,$ } 
		& \multirow{2}{*}{$\, \hs_\lambda(\ma,\mb) \,$}    & \multirow{2}{*}{$\, \os_\lambda(\ma,\mb) \,$ } \\[-5pt]
		&  &  &  &  &  &  \\
		\hline
		$\mu$     & 5.4  & 4.2  & 6.0  & 4.1  & 5.5  & 5.8  \\
		$\mu_1$   & 23.2 & 22.9 & 27.3 & 18.4 & 21.0 & 26.8 \\
		$\mu_2$   & 84.6 & 86.5 & 90.6 & 83.4 & 83.6 & 91.5 \\
		$\mu_3$   & 100  & 100  & 100  & 100  & 100  & 100  \\
		$\Sigma_1$& 7.7  & 7.0  & 7.2  & 7.8  & 7.0  & 4.9  \\
		$\Sigma_2$& 15.7 & 16.9 & 14.8 & 13.6 & 12.0 & 6.5  \\
		$\Sigma_3$& 100  & 100  & 100  & 94.8 & 100  & 94.8 \\ 
		\hline
		$avg.\,time$ & 0.80 & 0.63 & 0,73 & 0.79 & 0.66 & 0,75 \\
		\hline		
	\end{tabular}}
	\caption{Empirical power of the tests, measured in \% for simulated experiments on $k=10$ clusters, sample size $n = m = 1000$ and $\lambda = 0.5$, $0.75$, $1$, $5$, $10$ and $50$. \textit{avg. time} is measured in seconds.}
	\label{fig: Table2}
\end{table}

For the experiment of Table \ref{fig: Table2}, the computational cost decreases
as the value of $\lambda$ grows. For instance, for $\lambda=0.5$, the average time of computing a single $p$-value for the tests based on $\, S_\lambda(\ma,\mb) \,$, $\, \os_\lambda(\ma,\mb) \,$, and $\, \hs_\lambda(\ma,\mb) \,$ are $0.959$, $2.2940$ and $0.799$ seconds, respectively, while, for $\lambda=5$, the corresponding
computation times are $ 0.7701$, $0.7329$ and $0.6414$ seconds, respectively. 
When $\lambda$ goes above $5.0$, the reduction in computational cost becomes
relatively small.

\subsubsection{Variations on the samples' sizes.}

To evaluate the effect of sample size on the statistics considered, we work again in the setting of the MVG distribution, with $\lambda=1$, $k=10$, and let the sample sizes take the values $n=m=100$, $n=m=200$, $n=m=500$, $n=m=1000$, $n=m=2000$ and $n=m=5000$. Table \ref{fig: Table3} displays the power figures obtained in this simulation. In the case of alternatives with difference in mean, the power numbers show important increments even when we move from sample sizes 100 to 200, and continue to improve noticeably as the sample sizes increase through all values considered in the table. For the alternatives with difference in covariance matrices, the sample sizes required to achieve important increments in power are larger, a result that agrees with general belief, in the sense that changes in covariance are harder to detect.

\begin{table}[]
	\centering
	\resizebox{0.9\textwidth}{!}{%
	\begin{tabular}{|l|ccc|ccc|}
		\hline
		& \multicolumn{3}{c|}{$n,m=100$} & \multicolumn{3}{c|}{$n,m=200$}  \\ \cline{2-7} 
		& \multirow{2}{*}{ $\, S_\lambda(\ma,\mb) \,$ }  & \multirow{2}{*}{$\, \hs_\lambda(\ma,\mb) \,$}
		& \multirow{2}{*}{$\, \os_\lambda(\ma,\mb) \,$ } & \multirow{2}{*}{ $\, S_\lambda(\ma,\mb) \,$ } 
		& \multirow{2}{*}{$\, \hs_\lambda(\ma,\mb) \,$}    & \multirow{2}{*}{$\, \os_\lambda(\ma,\mb) \,$ } \\[-5pt]
		&  &  &  &  &  &  \\
		\hline
		$\mu$     & 6.8  & 4.5  & 5.2  & 5.0  & 5.4  & 4.0  \\
		$\mu_1$   & 8.9  & 7.2  & 7.0  & 8.9  & 9.1  & 7.3  \\
		$\mu_2$   & 17.2 & 18.5 & 16.2 & 34.3 & 28.2 & 34.3 \\
		$\mu_3$   & 100  & 100  & 100  & 100  & 100  & 100  \\
		$\Sigma_1$& 4.4  & 5.9  & 6.2  & 6.5  & 5.8  & 5.5  \\
		$\Sigma_2$& 6.2  & 7.8  & 6.9  & 7.7  & 8.7  & 5.9  \\
		$\Sigma_3$& 61.7 & 55.2 & 51.4 & 92.3 & 87.7 & 84.5 \\ 
		\hline
		$avg.\,time$ & 0.60 & 0.51 & 0.89 & 0.67 & 0.52 & 0.92 \\
		\hline\hline
		& \multicolumn{3}{c|}{$n,m=500$} & \multicolumn{3}{c|}{$n,m=1000$}  \\ \cline{2-7} 
		& \multirow{2}{*}{ $\, S_\lambda(\ma,\mb) \,$ }  & \multirow{2}{*}{$\, \hs_\lambda(\ma,\mb) \,$}
		& \multirow{2}{*}{$\, \os_\lambda(\ma,\mb) \,$ } & \multirow{2}{*}{ $\, S_\lambda(\ma,\mb) \,$ } 
		& \multirow{2}{*}{$\, \hs_\lambda(\ma,\mb) \,$}    & \multirow{2}{*}{$\, \os_\lambda(\ma,\mb) \,$ } \\[-5pt]
		&  &  &  &  &  &  \\
		\hline
		$\mu$     & 4.0  & 5.2  & 5.0  & 5.2  & 4.3  & 5.6  \\
		$\mu_1$   & 16.5 & 14.9 & 15.0 & 33.1 & 29.3 & 30.0 \\
		$\mu_2$   & 73.3 & 70.3 & 70.7 & 98.6 & 95.3 & 94.7 \\
		$\mu_3$   & 100  & 100  & 100  & 100  & 100  & 100  \\
		$\Sigma_1$& 7.3  & 6.3  & 6.5  & 8.8  & 8.5 & 7.5  \\
		$\Sigma_2$& 11.9 & 10.0 & 9.2  & 18.3 & 14.0 & 11.7 \\
		$\Sigma_3$& 100  & 99.9 & 99.7 & 100  & 100  & 100  \\
		\hline
		$avg.\,time$ & 0.71 & 0.69 & 0.84 & 1.08 & 0.89 & 1.17 \\
		\hline\hline
		& \multicolumn{3}{c|}{$n,m=2000$} & \multicolumn{3}{c|}{$n,m=5000$}  \\ \cline{2-7} 
		& \multirow{2}{*}{ $\, S_\lambda(\ma,\mb) \,$ }  & \multirow{2}{*}{$\, \hs_\lambda(\ma,\mb) \,$}
		& \multirow{2}{*}{$\, \os_\lambda(\ma,\mb) \,$ } & \multirow{2}{*}{ $\, S_\lambda(\ma,\mb) \,$ } 
		& \multirow{2}{*}{$\, \hs_\lambda(\ma,\mb) \,$}    & \multirow{2}{*}{$\, \os_\lambda(\ma,\mb) \,$ } \\[-5pt]
		&  &  &  &  &  &  \\
		\hline
		$\mu$     & 6.1  & 5.6  & 4.8  & 4.7  & 5.6  & 5.7  \\
		$\mu_1$   & 60.5 & 54.7 & 55.1 & 97.5 & 94.9 & 93.9 \\
		$\mu_2$   & 100  & 99.8 & 100  & 100  & 100  & 100  \\
		$\mu_3$   & 100  & 100  & 100  & 100  & 100  & 100  \\
		$\Sigma_1$& 10.9 & 8.9  & 7.4  & 22.9 & 19.9 & 17.4  \\
		$\Sigma_2$& 28.5 & 19.8 & 17.2 & 75.4 & 42.8 & 37.8 \\
		$\Sigma_3$& 100  & 100  & 100  & 100  & 100  & 100  \\
		\hline
		$avg.\,time$ & 1.47 & 1.24 & 1.53 & 2.81 & 2.57 & 2.99 \\
		\hline
	\end{tabular}}
	\caption{Empirical power of the tests, measured in \% for simulated experiments on $k=10$ clusters, $\lambda = 1$, and sample sizes $n = m =100$, $200$, $500$, $1000$, $2000$ and $5000$. \textit{avg. time} is measured in seconds. }
	\label{fig: Table3}
\end{table}

The computation time increases with sample size, as expected, but the savings
by using smaller sample sizes is not that significant, in the range of values
considered in our experiment. For instance, for $n,m= 100$, the average $p$-value computation times for $\, S_\lambda(\ma,\mb) \,$, $\, \os_\lambda(\ma,\mb) \,$, and $\, \hs_\lambda(\ma,\mb) \,$ are $0.6018$, $0.9253$ and $0.5120$, respectively,
which do not represent important savings respect to the times reported above
for $n,m= 500$. This suggests that for small sample sizes, the ``fixed costs'' 
associated with the setting up of the statistic is an important part of the
computational cost. On the same line, a relatively large increase of the sample sizes, respect to those considered in Table \ref{fig: Table1}, does not 
significantly affect the computation time. For instance, for $n=m=5000$, the average computation time of the $p$-values of the $\, S_\lambda(\ma,\mb) \,$, $\, \os_\lambda(\ma,\mb) \,$, and $\, \hs_\lambda(\ma,\mb) \,$ statistic are $1.4725$, $1.5367$ and $1.2413$.

\subsubsection{Variations on the dimension of the dataset.}

To evaluate the effect of changing the dimension of the datasets on the statistics studied, we use the bounded Burr distributions, since these depend only on the 
unidimensional parameter $\alpha$, making the generation process simpler. 

\begin{table}[]
	\centering
	\resizebox{0.9\textwidth}{!}{%
	\begin{tabular}{|c|ccc|ccc|}
		\hline
		& \multicolumn{3}{c|}{$dim = 5$} & \multicolumn{3}{c|}{$dim = 10$}  \\ \cline{2-7} 
		& \multirow{2}{*}{ $\, S_\lambda(\ma,\mb) \,$ }  & \multirow{2}{*}{$\, \hs_\lambda(\ma,\mb) \,$}
		& \multirow{2}{*}{$\, \os_\lambda(\ma,\mb) \,$ } & \multirow{2}{*}{ $\, S_\lambda(\ma,\mb) \,$ } 
		& \multirow{2}{*}{$\, \hs_\lambda(\ma,\mb) \,$}    & \multirow{2}{*}{$\, \os_\lambda(\ma,\mb) \,$ } \\[-5pt]
		&  &  &  &  &  &  \\
		\hline
		$\alpha_1$ & 5.5  & 5.3  & 4.9  & 4.3  & 4.3  & 4.3  \\
		$\alpha_2$ & 11.7 & 10.0 & 13.0 & 15.8 & 9.3  & 14.9 \\
		$\alpha_3$ & 30.9 & 24.1 & 33.6 & 45.6 & 39.8 & 48.9 \\
		$\alpha_4$ & 57.9 & 53.5 & 61.9 & 79.4 & 73.4 & 83.3 \\
		$\alpha_5$ & 78.0 & 78.3 & 81.8 & 93.9 & 94.1 & 96.5 \\
		$\alpha_6$ & 12.2 & 12.2 & 14.2 & 19.2 & 14.1 & 19.4 \\
		$\alpha_7$ & 87.1 & 74.4 & 80.1 & 92.1 & 83.0 & 87.5 \\
		\hline
		$avg.\,time$ & 0.79 & 0.63 & 0.76 & 0.84  & 0.79 & 0.91 \\
		\hline\hline
		& \multicolumn{3}{c|}{$dim = 20$} & \multicolumn{3}{c|}{$dim = 50$}  \\ \cline{2-7} 
		& \multirow{2}{*}{ $\, S_\lambda(\ma,\mb) \,$ }  & \multirow{2}{*}{$\, \hs_\lambda(\ma,\mb) \,$}
		& \multirow{2}{*}{$\, \os_\lambda(\ma,\mb) \,$ } & \multirow{2}{*}{ $\, S_\lambda(\ma,\mb) \,$ } 
		& \multirow{2}{*}{$\, \hs_\lambda(\ma,\mb) \,$}   & \multirow{2}{*}{$\, \os_\lambda(\ma,\mb) \,$ } \\[-5pt]
		&  &  &  &  &  &  \\
		\hline
		$\alpha_1$ & 4.5  & 4.3  & 5.8  & 4.8  & 5.8  & 5.8  \\
		$\alpha_2$ & 14.8 & 13.5 & 13.4 & 17.1 & 15.5 & 14.5 \\
		$\alpha_3$ & 55.1 & 45.2 & 53.3 & 63.3 & 56.0 & 56.7 \\
		$\alpha_4$ & 88.2 & 84.4 & 90.5 & 92.9 & 91.5 & 93.0 \\
		$\alpha_5$ & 98.7 & 98.2 & 99.4 & 99.4 & 99.4 & 99.6 \\
		$\alpha_6$ & 21.4 & 18.9 & 17.8 & 20.4 & 22.6 & 21.2 \\
		$\alpha_7$ & 94.0 & 87.9 & 91.0 & 94.2 & 91.8 & 91.8 \\		
		\hline
		$avg.\,time$ & 0.96 & 0.91 & 1.01 & 1.35 & 1.30 & 1.39 \\
		\hline\hline
		& \multicolumn{3}{c|}{$dim = 100$} & \multicolumn{3}{c|}{$dim = 200$}  \\ \cline{2-7} 
		& \multirow{2}{*}{ $\, S_\lambda(\ma,\mb) \,$ }  & \multirow{2}{*}{$\, \hs_\lambda(\ma,\mb) \,$}
		& \multirow{2}{*}{$\, \os_\lambda(\ma,\mb) \,$ } & \multirow{2}{*}{ $\, S_\lambda(\ma,\mb) \,$ } 
		& \multirow{2}{*}{$\, \hs_\lambda(\ma,\mb) \,$}   & \multirow{2}{*}{$\, \os_\lambda(\ma,\mb) \,$ } \\[-5pt]
		&  &  &  &  &  &  \\
		\hline
		$\alpha_1$ & 5.2  & 4.8  & 4.5  & 4.2  & 4.4  & 4.7  \\
		$\alpha_2$ & 14.4 & 16.0 & 15.6 & 16.1 & 18.0 & 16.9 \\
		$\alpha_3$ & 64.3 & 59.3 & 58.9 & 64.2 & 66.2 & 67.3 \\
		$\alpha_4$ & 94.9 & 95.0 & 94.0 & 95.0 & 94.9 & 93.9 \\
		$\alpha_5$ & 99.8 & 99.6 & 99.6 & 100  & 100  & 99.7 \\
		$\alpha_6$ & 24.6 & 22.6 & 21.2 & 23.2 & 24.7 & 21.7 \\
		$\alpha_7$ & 95.1 & 92.7 & 92.6 & 94.0 & 93.7 & 95.4 \\
		\hline
		$avg.\,time$ & 2.75 & 2.39 & 2.42 & 3.87 & 3.86 & 4.17 \\
		\hline
	\end{tabular}}
	\caption{Empirical power of the tests, measured in \% for simulated experiments on $k=10$ clusters, $\lambda = 1$, and sample sizes $n = m =500$ and $d=5$, $10$, $20 $, $50$, $100$ and $200$. \textit{avg. time} is measured in seconds.}
	\label{fig: Table4}
\end{table}

Table \ref{fig: Table4} shows the changes in power and average computation time, when 
$\alpha$ varies over the same set of values used in Table \ref{fig: Table1}, while
the dimension takes values $d=5, 10, 20, 50, 100$ and $200$.  As expected, as the 
dimension of the samples increase, the power of the statistics increases considerably, over the range of dimensions included. An exception to this tendency is the behavior of $S_\lambda(\ma,\mb)$, whose power does not improve when dimension goes beyond $d=50$. Interestingly, plots of 
the average computation time suggest that computation times of the statistics
considered grow sub-linearly with dimension, making
these statistics applicable in large dimensions.

\subsubsection{Variations on the number of clusters.}

An important parameter in our procedure is the number $k$, of clusters, in the
$k$-means clustering procedure applied to the pooled sample. Increasing $k$
should improve the ability of the procedure to distinguish between different distributions, but at the same time, increasing $k$ directly increases the dimension
of the optimization problem to be solved in the Sinkhorn algorithms. To evaluate
the effect of variations in $k$ we consider the Dirichlet distribution described
above, in dimension $d=5$, with the same parameters, for the null and alternatives, listed in Table \ref{fig: Table1}, letting $k$ take the values $5$, $10$, $15$, 
$20$, $25$ and $30$.  Table \ref{fig: Table5} shows the corresponding estimated
power numbers and average computation times (for one $p$-value) for this 
experiment.

\begin{table}[]
	\centering
	\resizebox{0.9\textwidth}{!}{%
	\begin{tabular}{|c|ccc|ccc|}
		\hline
		& \multicolumn{3}{c|}{$k = 5$} & \multicolumn{3}{c|}{$k = 10$}  \\ \cline{2-7} 
		& \multirow{2}{*}{ $\, S_\lambda(\ma,\mb) \,$ }  & \multirow{2}{*}{$\, \hs_\lambda(\ma,\mb) \,$}
		& \multirow{2}{*}{$\, \os_\lambda(\ma,\mb) \,$ } & \multirow{2}{*}{ $\, S_\lambda(\ma,\mb) \,$ } 
		& \multirow{2}{*}{$\, \hs_\lambda(\ma,\mb) \,$}   & \multirow{2}{*}{$\, \os_\lambda(\ma,\mb) \,$ } \\[-5pt]
		&  &  &  &  &  &  \\
		\hline
		$\beta_1$ & 4.4  & 4.8  & 4.8  & 4.5  & 5.7  & 4.5  \\
		$\beta_2$ & 10.0 & 98.7 & 23.7 & 83.9 & 91.4 & 33.9 \\
		$\beta_3$ & 54.8 & 100  & 100  & 100  & 100  & 100  \\
		$\beta_4$ & 10.5 & 13.7 & 13.1 & 9.5  & 12.7 & 14.4 \\
		$\beta_5$ & 100  & 100  & 100  & 100  & 100  & 100  \\
		$\beta_6$ & 32.9 & 100  & 89.5 & 100  & 100  & 98.7 \\ 
		\hline
		$avg.\,time$ & 0.55 & 0.53 & 0.61 & 0.89 & 0.67 & 0.77 \\
		\hline\hline
		& \multicolumn{3}{c|}{$k = 15$} & \multicolumn{3}{c|}{$k = 20$}  \\ \cline{2-7} 		
		& \multirow{2}{*}{ $\, S_\lambda(\ma,\mb) \,$ }  & \multirow{2}{*}{$\, \hs_\lambda(\ma,\mb) \,$}
		& \multirow{2}{*}{$\, \os_\lambda(\ma,\mb) \,$ } & \multirow{2}{*}{ $\, S_\lambda(\ma,\mb) \,$ } 
		& \multirow{2}{*}{$\, \hs_\lambda(\ma,\mb) \,$}   & \multirow{2}{*}{$\, \os_\lambda(\ma,\mb) \,$ } \\[-5pt]
		&  &  &  &  &  &  \\
		\hline
		$\beta_1$ & 5.0  & 5.2  & 6.1  & 4.9  & 4.8  & 4.7  \\
		$\beta_2$ & 77.1 & 92.4 & 44.8 & 86.1 & 86.1 & 46.5 \\
		$\beta_3$ & 100  & 100  & 100  & 100  & 100  & 100  \\
		$\beta_4$ & 6.7  & 13.6 & 13.8 & 9.0  & 13.8 & 14.1 \\
		$\beta_5$ & 100  & 100  & 100  & 100  & 100  & 100  \\
		$\beta_6$ & 100  & 100  & 99.8 & 100  & 100  & 99.9 \\ 
		\hline
		$avg.\,time$ & 1.03 & 1.00 & 1.15 & 1.32 & 1.27 & 1.51 \\
		\hline\hline
		& \multicolumn{3}{c|}{$k = 25$} & \multicolumn{3}{c|}{$k = 30$}  \\ \cline{2-7} 
		& \multirow{2}{*}{ $\, S_\lambda(\ma,\mb) \,$ }  & \multirow{2}{*}{$\, \hs_\lambda(\ma,\mb) \,$}
		& \multirow{2}{*}{$\, \os_\lambda(\ma,\mb) \,$ } & \multirow{2}{*}{ $\, S_\lambda(\ma,\mb) \,$ } 
		& \multirow{2}{*}{$\, \hs_\lambda(\ma,\mb) \,$}   & \multirow{2}{*}{$\, \os_\lambda(\ma,\mb) \,$ } \\[-5pt]
		&  &  &  &  &  &  \\
		\hline
		$\beta_1$ & 4.6  & 5.4  & 5.1  & 4.4  & 5.3  & 5.3  \\
		$\beta_2$ & 89.3 & 79.5 & 44.7 & 89.6 & 78.9 & 48.6 \\
		$\beta_3$ & 100  & 100  & 100  & 100  & 100  & 100  \\
		$\beta_4$ & 5.8  & 14.8 & 14.2 & 7.4  & 12.5 & 13.4 \\
		$\beta_5$ & 100  & 100  & 100  & 100  & 100  & 100  \\
		$\beta_6$ & 100  & 100  & 100  & 100  & 100  & 100  \\ 
		\hline
		$avg.\,time$ & 1.66 & 1.52 & 1.93 & 2.01 & 1.91 & 2.58  \\
		\hline
	\end{tabular}}
	\caption{Empirical power of the tests, measured in \% for simulated experiments with parameter $\lambda=1$ and $k = 5$, $10$, $15$, $20$, $25$ and $30$ clusters. \textit{avg. time} is measured in seconds.}
	\label{fig: Table5}
\end{table}

In Table \ref{fig: Table5} it is evident that, in general (with the 
exception of the most difficult alternatives), the power of the statistics improves
very noticeably when $k$ goes from $5$ to $10$. Another noticeable increase
in power occurs when going from $k=10$ to $15$, although this second increment
is less important. But going beyond, $k=15$ does not produce important changes
in the power of the statistics, for this type of data. Thus, for this parameter
our recommendation would be to stay with $k=10$ or $k=15$, at most, although
this conclusion could be affected by changes in the data dimension and
sample sizes. Again, as
happened in the case of changes in data dimension, plotting the average 
computation times against the values of $k$, a nearly linear tendency is
observed, suggesting that the computational cost is linear in the number
of clusters employed.

{\subsection{Comments on the Power Comparisons}}
The three different possible implementations considered of permutation tests based
on the Sinkhorn divergence compare favourably, in terms of power, with the classical
and asymptotically efficient test of Schilling, and compare favourably, in terms of 
power and computational cost, with the test based on the classical optimal 
transport statistic (Wasserstein). All three exhibit good power against the varied 
set of distributional scenarios considered. The computational cost for these tests 
allows their practical application for sample sizes in the order of 5000 and 
dimensions up to $d=50$, at least. 

For the scenarios analysed, based on considerations of statistical power
and computational cost, an intermediate value of $\lambda$, between 1 and 10,
seems to be a good choice for this parameter, while, the same considerations 
would suggest that $k=10$ is a good selection for the number of clusters in
the $k$-means clustering procedure. 

Overall, and ignoring some particular cases, there are not very important 
differences, in terms of power, between the different versions of the Sinkhorn 
statistic considered. Still, our preference among them, would go with 
$\hs_{\lambda}(\ma,\mb)$, for its overall statistical power, lower 
computational cost and for having more results available in the literature, 
including results on the approximation of the solution to the classical OT 
problem by $\hs_{\lambda}(\ma,\mb)$, for an appropriate choice of $\lambda$ 
(see \cite{awr17}).

	\subsection{Evaluating the Convergence to Normality}
	The purpose of this subsection is to evaluate, in an example, the convergence
	stated in Theorem \ref{clt1} for $\hs_\lambda(\hat{\ma}_N,\hat{\mb}_N)$ under the null hypothesis and discuss the practical approximation of the limiting
	parameters appearing in (\ref{limit1}). 
	
\begin{figure}[]
	\centering
	\includegraphics[scale=0.37]{"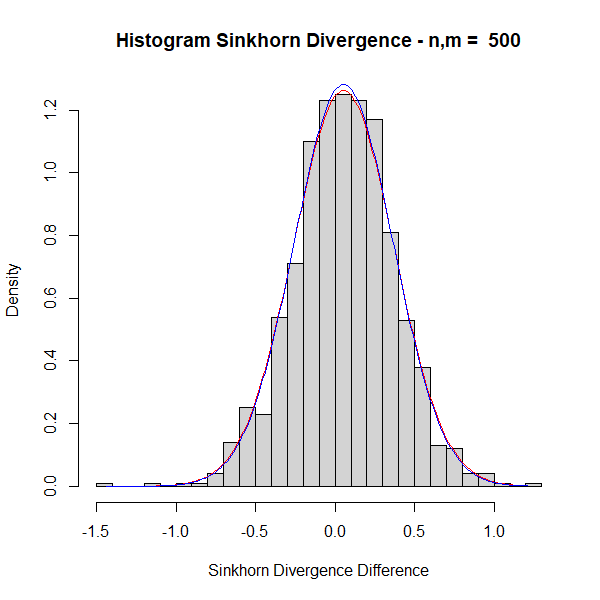"} 
	\includegraphics[scale=0.37]{"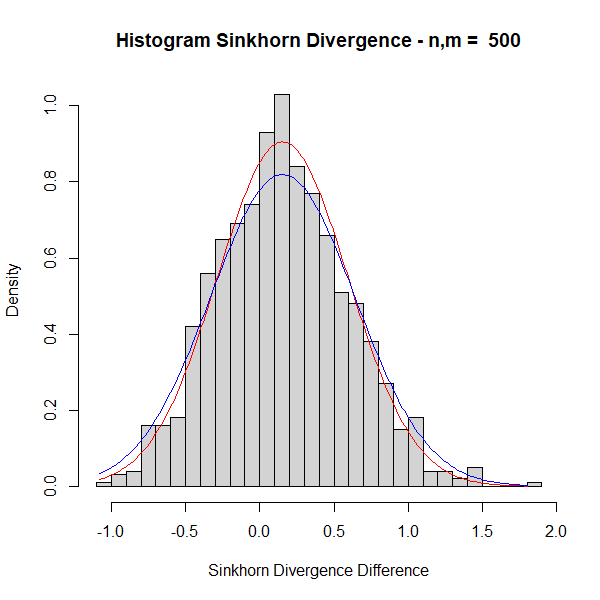"} 
	\includegraphics[scale=0.37]{"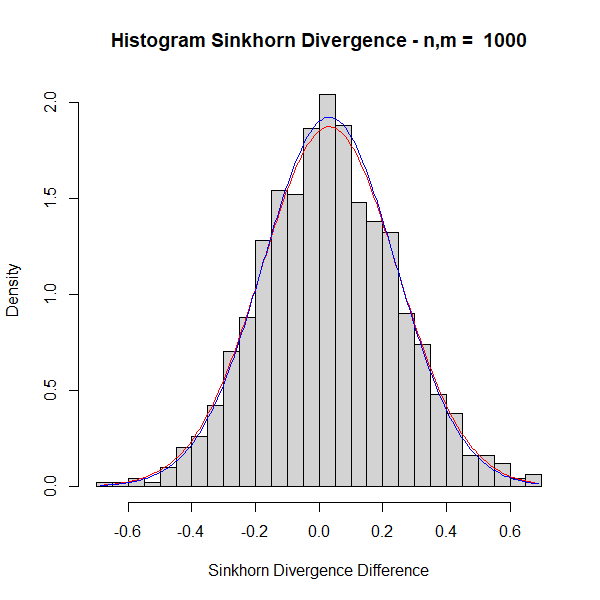"} 
	\includegraphics[scale=0.37]{"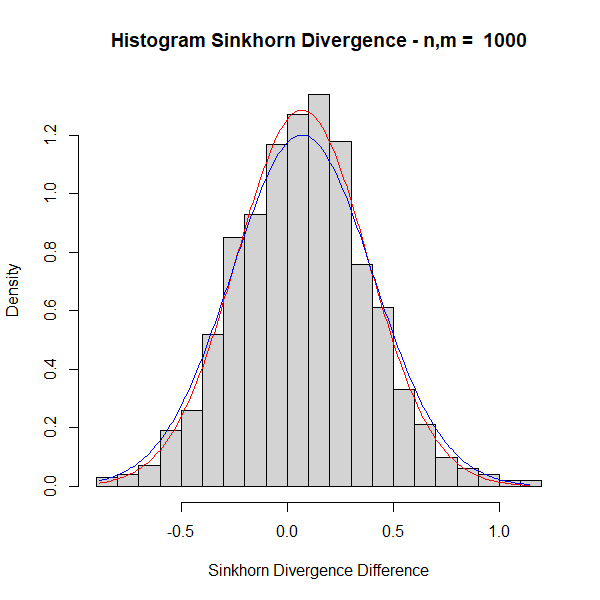"} 
	\includegraphics[scale=0.37]{"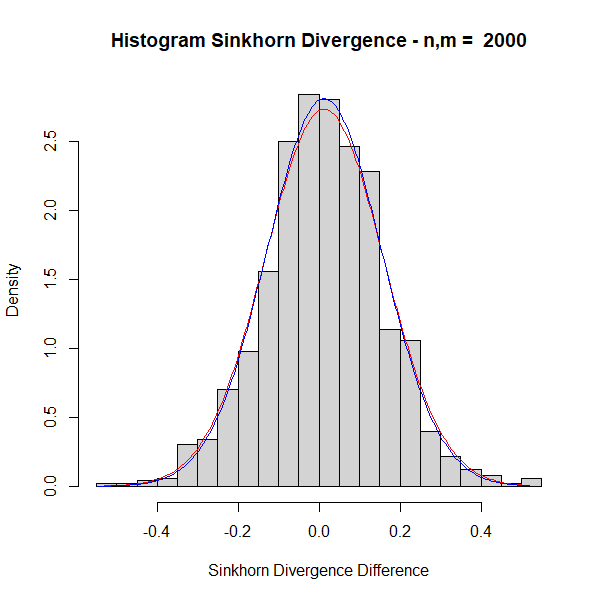"}
	\includegraphics[scale=0.37]{"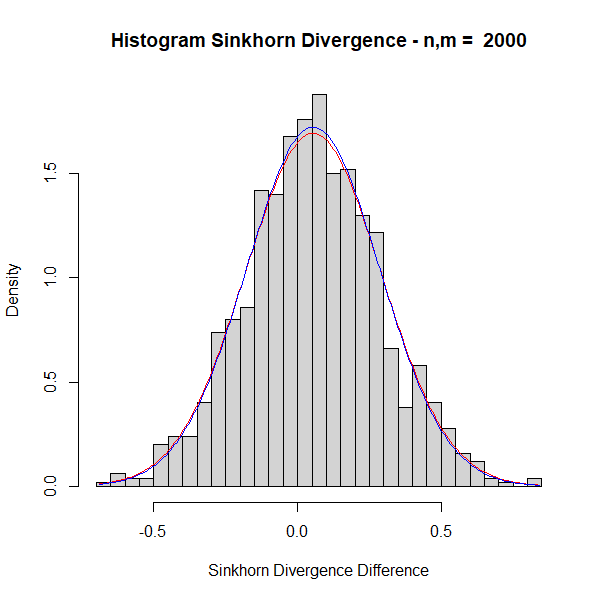"} 
	\caption{Histograms of {\small$\sqrt{N_e}(\hs_\lambda(\hat{\ma}_N,\hat{\mb}_N)-\hs_\lambda({\ma}^*_N,{\ma}^*_N))$} in the case of MVG samples under the null hypothesis. On the left column 
	we find the histograms obtained with $k=10$ and on the right column $k=15$.} 
	\label{fig:SinkDivHist}
\end{figure}

	We work again with i.i.d. samples, $X$ and $Y$, from the multivariate standard Gaussian distribution. Based on the results of the previous subsections, we
	use $\lambda=10$,  with number of 	clusters $k=10$ and $k=15$, and let the sample sizes take the values $n=m=500$, $n=m=1000$ and $n=m=2000$. For each pair of samples, we proceed as follows: First, the $k$-means clusters for the combined sample are computed, as in the previous experiments. Then, the population (limiting) probabilities for these cells, that form the vector ${\ma}^*_N$,  
	are estimated from a considerably larger i.i.d. sample, $X_0$, of sample size
	    $250000$, and these approximated cell probabilities
	are used to compute the matrix $C$ in (\ref{cova}). Next, to compute 
	the gradient vector $J_I$ in (\ref{limit1}), for each $i\leq k-1$, the
	vector $\ma_{N,\,\varepsilon}^*(i)$ is obtained from $\ma^*_N$ by substracting a small amount $\varepsilon$ from its $i$-th component and adding the same amount to its $k$-th component (so that the
	vector remains a probability vector) and the $i$-th	partial derivative, with 
	respect to the coordinates of $\ma^*_N$ is
	computed as
	\begin{equation}
		\begin{split}
			&\dfrac{\,\hs_\lambda({\ma}^*_N,{\ma}^*_N)-\hs_\lambda(\ma_{N,\,\varepsilon}^*(i),\ma_{N}^*)\,}{\varepsilon},\quad \text{for } 1\leq i \leq k-1 \\
			\text{ and }\quad  & \dfrac{\,\hs_\lambda({\ma}^*_N,{\ma}^*_N) - \hs_\lambda({\ma}^*_N,\ma_{N,\,\varepsilon}^*(i))\,}{\varepsilon},\quad \text{for } k+1\leq i \leq 2k-1.
		\end{split}
	\end{equation}
	In our estimation, $\varepsilon$ is set to $10^{-7}$. With $\ma_N^*$, $C$ and
	$J_I$, both the left and right side in statement (\ref{limit1}) are approximated. For each choice of parameters, these calculations are repeated $1000$ times and the results are reported in Figure \ref{fig:SinkDivHist}, where
	the histograms correspond to the variable $\nu_N=\sqrt{N_{\mbox{e}}}(\hs_\lambda(\hat{\ma}_N,\hat{\mb}_N)-\hs_\lambda({\ma}^*_N,{\ma}^*_N))$, the blue line represents the normal density with the sample mean and variance of the calculated $\nu_N$ values, and the red line represents the average Gaussian density predicted by Theorem \ref{clt1}
	(mean and variance averaged over the 1000 realizations).
It appears that in all cases the distribution of the statistic $\nu_N$ is close
to normality (the samples of 1000 $\nu_N$ values pass the Shapiro normality test) and close to the distribution predicted by the theorem. The agreement between the red and blue curves improves with sample size, being 
very good for $n=m=2000$.

Practical use of the convergence exhibited in this computational evaluation would be limited to the case
where, besides the samples being used for the test statistic calculation, there exists a larger sample from
the same $X$ and $Y$ distributions, to be used for the estimation of the vector $\ma_N^*$, the matrix $C$ and
the partial derivatives in $J_I$, a possibility that might actually occur in times of big data. In this case,
one could use the estimated Gaussian distribution of the statistic for the significance evaluation,
instead of the permutation test.



\begin{thebibliography}{99} 

\bibitem{awr17} Altschuler, J., Weed, J. and Rigollet, P. (2017). Near-linear time \linebreak approximation algorithms for optimal transport via Sinkhorn iteration. In \linebreak{\it Advances in Neural Information Processing Systems}, 30 (NIPS 2017), 1961-1971.

\bibitem{bigot} Bigot, J., Cazelles, E. and Papadakis, N. (2017) Central limit theorems for Sinkhorn divergence between probability distributions on finite spaces and statistical applications. arXiv preprint arXiv:1711.08947.

\bibitem{Chen} Chen, H. and Friedman, J. H. (2017) A New Graph-Based Two-Sample Test for Multivariate and Object Data. {\it Journal of the American Statistical Association}, 112:517, 397-409.


\bibitem{cuturi} Cuturi, M. (2013) Sinkhorn distances: Lightspeed computation of optimal transport. In {\it Proceedings of Advances in Neural Information Processing Systems}, pp. 2292–2300.

\bibitem{cuturi2014} Cuturi, M.and Doucet, A. (2014) Fast Computation of Wasserstein 
Barycenters. In {\it Proceedings of the 31st International Conference on \linebreak Machine Learning, PMLR} 32(2) pp. 685-693.

\bibitem{Del Barrio} Del Barrio, E., Cuesta-Albertos, J. A., Matr\'an, C. and Rodr{\ai}guez-Rodr{\ai}guez, J. M. (1999) Tests of goodness of fit based on the 
\linebreak L$^2$-Wasserstein distance. 
{\it Annals of Statistics}, {\bf 27}, No. 4, pp. 1230–1239. 

\bibitem{Del Barrio L} Del Barrio, E. and Loubes, J. M. (2019) Central limit theorems for \linebreak empirical transportation in general dimension.  {\it The Annals of 
Probability}, {\bf 47}, No. 2, pp. 926–951. 

\bibitem{Dudley1999} Dudley, R. M. (1999) {\em Uniform Central Limit Theorems}.
Cambridge \linebreak University Press, Cambridge.

\bibitem{fr79} Friedman, J. H. and Rafsky, L. C. (1979) Multivariate Generalizations 
of the Wald-Wolfowitz and Smirnov Two-Sample Tests {\em The Annals of Statistics}, 
{\bf 7}, No. 4, 697-717.

\bibitem{frogner} Frogner, C., Zhang, C., Mobahi, H., Araya-Polo, M. and Poggio, T. 
(2015). Learning with a Wasserstein Loss. In {\it Proceedings of Advances in Neural 
Information Processing Systems}, NIPS 2015.

\bibitem{genevay19} Genevay, A., Chizat, L., Bach, F., Cuturi, M. and Peyr\'e, G.
(2019) Sample conplexity of Sinkhorn divergences. In {\it Proceedings of the 22nd
International Conference on Artificial Intelligence and Statistics} (AISTATS).

\bibitem{good05} Good, P. (2005) {\em Permutation, Parametric and Bootstrap Tests of Hypothesis}. Springer, New York.

\bibitem{johnson} Johnson, M. E. (1987)
 {\em Multivariate Statistical Simulation: A guide to selecting and generating continuous multivariate distributions}. John Wiley \& Sons, New York.
 
\bibitem{kbj} Kotz, S., Balakrishnan, N. and Johnson, N. L. (2000). {\em Continuous Multivariate Distributions. Volume 1: Models and Applications}. John Wiley \& Sons, New York.

\bibitem{mena} Mena, G. and Niles-Weed, J. (2019) Statistical bounds for entropic
optimal transport: sample complexity and the central limit theorem. In {\it 
Advances in Neural Information Processing Systems}, {\bf 32}, NeurIPS 2019.

\bibitem{NilesWeedRigollet} Niles-Weed, J., and Rigollet, P. (2021),
Estimation of Wasserstein distances in the Spiked Transport Model,
\textit{Bernoulli}, to appear.

\bibitem{pollard82} Pollard, D. (1982) A Central Limit Theorem for k-means
clustering. {\it Annals of Probability}, {\bf 10}(4), pp. 919–926.

\bibitem{rokach05} Rokach, L. and Maimon, O. (2005) Clustering Methods. Chapter
15 in {\it Data Mining and Knowledge Discovery Handbook. O. Maimon and L. Rokach, editors}, pp. 321-352. Springer, Boston, MA. 

\bibitem{schilling} Schilling, M. F. (1986) Two-Sample Tests Based on Nearest Neighbors. {\it Journal of the American Statistical Association}, {\bf 81}, No. 395, pp. 799-806. 

\bibitem{sommerfeld} Sommerfeld, M. and Munk, A. (2018) Inference for empirical Wasserstein distances on finite spaces. {\it Journal of the Royal Statistical Society: Series B} (Statistical Methodology), {\bf 80}(1): 219–238.

\bibitem{van der vaart} van der Vaart, A. W. (1998) {\em Asymptotic Statistics}.
Cambridge University Press, Cambridge.

\bibitem{villani} Villani, C.  (2009) {\em Optimal Transport. Old and New}.
Volume 338 in series Grundlehren der matematischen Wissenschaften. Springer.
Berlin, Heidelberg.

\end{thebibliography}
\end{document}